\documentclass[11pt]{article}
\usepackage{amsfonts}
\usepackage{amscd,amssymb,stmaryrd}
\usepackage{graphicx}
\usepackage{color,epsfig,comment}

\def\no{\noindent}  
\def \Vh0{\stackrel{\circ}{V}_h}

\def\l{\label}    

\def\b{\beta}  \def\a{\alpha} 
\def\m{\mbox} \def\t{\times}
\def\Box{\sharp}

\newcommand{\stl}
{\mathrel{\raise2pt\hbox{${\mathop<\limits_{\raise1pt\hbox{\mbox{$\sim$}}}}$}}}

\newcommand{\stg}
{\mathrel{\raise2pt\hbox{${\mathop>\limits_{\raise1pt\hbox{\mbox{$\sim$}}}}$}}}

\newcommand{\ste}
{\mathrel{\raise2pt\hbox{${\mathop=\limits_{\raise1pt\hbox{\mbox{$\sim$}}}}$}}}

\newtheorem{lemma}{Lemma}[section]

\newtheorem{theorem}{Theorem}[section]

\newtheorem{remark}{Remark}[section]
\def\bd{\begin{description}} \def\ed{\end{description}}

\def \no{\noindent}

\def \ee{\begin{equation}}
\def \e{\end{equation}}
\def\beq{\begin{eqnarray}}
\def\eq{\end{eqnarray}}
\def\beqx{\begin{eqnarray*}}
\def\eqx{\end{eqnarray*}}
\def\l{\label}  \def\no{\noindent}

\def\omega{\alpha}

\def\ff{ {\mbox{\sc f}}}
\def\E{{\mbox{\sc e}}}
\def\vv{\mbox{\sc v}}
\def\o{\Omega}
\def\c{{\bf curl}}

\def\v{{\bf v}}

\def\x{{\bf x}}
\def\ti{\times}
\def\la{{\bf \lambda}}

\def\w{{\bf w}}
\def\n{{\bf n}}

\def\0{{\bf 0}}
\def\a{{\bf a}}
\def\b{{\bf b}}
\def\t{{\bf t}}
\def\r{{\bf r}}

\def\T{{\cal T}_{h}}
\def\12{{1\over 2}}

\def\rig{\rightarrow}

\def\chi{{\cal X}}

\def\cu{{\bf curl}}

\begin{document}

\begin{center}
{\large\bf Convergence of the HX Preconditioner for Maxwell's Equations with Jump Coefficients (i): Extensions of the Regular Decomposition}
\end{center}

\bigskip
\centerline{Qiya ~Hu\footnote{1. LSEC, ICMSEC, Academy of Mathematics and Systems Science, Chinese Academy of Sciences, Beijing 100190, China; 2. School of Mathematical Sciences, University of Chinese Academy
of Sciences, Beijing 100049, China (hqy@lsec.cc.ac.cn). This author was supported by the
Natural Science Foundation of China G11571352.}}
\bigskip

\medskip
\begin{abstract}
This paper is the first of the two serial articles, aiming to prove the convergence of the HX preconditioner for Maxwell's equations with jump coefficients,
originally proposed by Hiptmair and Xu \cite{HX}. In this paper, we establish several extensions of the discrete regular
decomposition for edge finite element functions defined in three-dimensional domains. The functions defined by the new discrete regular decompositions can inherit
zero degrees of freedom of the considered edge finite element function on some faces and edges of polyhedral domains as well as of some non-Lipschitz domains, and possess nearly
optimal stability with only a $logarithmic$ factor. These regular decompositions will be used to prove the convergence of the HX preconditioner for Maxwell's
equations with jump coefficients in a subsequent paper \cite{Hu2-2017}.
\end{abstract}

\medskip

{\bf AMS(MOS) subject classification}. 65N30, 65N55

\tableofcontents

\section{Introduction}\l{sec:introduction}
\setcounter{equation}{0}

It is known that any vector-valued function in $H(\c)$ space can be decomposed into the sum of a $H^1$ vector-valued function and the gradient
of a $H^1$ scalar-valued function (refer to \cite{Gir1986} and \cite{Regul}), and the decomposition is stable with respect to the standard norms. If the sum is orthogonal with respect to the $L^2$ inner product,
the decomposition is called the Helmholtz decomposition, otherwise, it is called a regular decomposition (which is not unique). The regular decompositions
are more popular than the Helmholtz decomposition because the regular decompositions are valid on the general Lipschitz domains, whereas the Helmholtz decomposition
holds only on smooth or convex domains (if $H^1$-regularity of the vector-valued function is not required, the Helmholtz decomposition can always be found for the general Lipschitz domains).
In the study of Maxwell's equations, the Helmholtz decompositions and regular decompositions can be used to transform a problem for a $H(\c)$ function into problems on two $H^1$ functions.

The Nedelec edge finite element method (refer to \cite{s8}) is a popular discretization method for Maxwell's equations.
Suitable decompositions for edge finite element functions, known as discrete (Helmholtz or regular) decompositions, play a key role
in the convergence analysis of preconditioners for Maxwell's equations (see, for example, \cite{hipt,HX,HuShuZou2013,hz,HZ2,PasciakZhao2002,s10,Tos}).
In this analysis, Maxwell's equations in a {\it homogeneous} medium are often assumed because the stability of the standard Helmholtz decomposition only holds for the usual norms.
However, nonhomogeneous media often occur in practical applications; therefore, some
weighted norms have to be introduced. A natural question is: whether the discrete decompositions are still stable with respect to the weighted norms?
This question was first positively answered in a study \cite{HuShuZou2013} by constructing a discrete {\it weighted orthogonal} Helmholtz decomposition
that was proved to be almost stable with respect to weighted norms.

For the numerical analysis of Maxwell's equations in a nonhomogeneous medium, a discrete {\it regular} decomposition needs to be developed for edge finite element functions defined in three-dimensional
domains such that the decomposition is stable with respect to weighted norms. This study is complex; therefore, it will be published in two parts.
Herein, some technical tools are developed to derive some extensions of the discrete regular decomposition for edge finite element functions.
The standard regular (and Helmholtz) decompositions possess a vital property: when the vector-valued function being considered has zero trace on the boundary of the underlying
domain, the functions defined by the decomposition also have zero trace on this boundary. We will construct a discrete {\it regular}  decomposition on a polyhedral domain
such that the inheritance of the zero trace on the boundary can be maintained when the boundary is replaced by a union of some local faces and edges of
the polyhedral domain. In particular, we also establish the corresponding decompositions for edge finite element function on some non-Lipschitz domains,
which are unions of two polyhedral domains intersecting at one edge or vertex.
We will show that the resulting regular decompositions possess stability estimates
with only one $logarithmic$ factor. These interesting results will be used in a subsequent article \cite{Hu2-2017} to develop a
discrete regular decomposition that is nearly stable with respect to weighted norms, and the convergence of the HX preconditioner for the case with
jump coefficients will be further proved.

The remainder of this paper is organized as follows. Section ~2 defines some edge finite element subspaces. Section 3 presents the main results. Section 4 discusses the regular decompositions that preserve zero-tangential traces on faces. Section 5 presents several discrete regular decompositions that preserve local zero-tangential components on edges and faces; moreover, this section proves the result for polyhedral domains as well. In Section 6, discrete regular decompositions on some non-Lipschitz domains are derived.

\section{Preliminaries}
\setcounter{equation}{0}

In this section we introduce some fundamental notions and notations.

\subsection{Sobolev spaces and norms}

For an open and connected bounded domain $G$ in $\mathbf
R^3$, let $H^1(G)$ be the standard Sobolev space consisting of the functions whose weak derivatives belong to $L^2(G)$.
Define the $H({\bf curl})$ space as follows
$$
H({\bf curl}; G)=\{\v\in L^2(G)^3; ~\cu\,\v\in
L^2(G)^3\}.$$
For the two spaces, we use natural norms
$$ \|\v\|_{1,G}=(|\v|^2_{1,G}+\|\v\|^2_{0,G})^{{1\over 2}},~~~\v\in (H^1(G))^3 $$
and
$$ \|\v\|_{\c,G}=(\|\c~\v\|^2_{0,G}+\|\v\|^2_{0,G})^{{1\over 2}},~~~\v\in H(\c;~G). $$

For a polyhedron $G$, let $\Gamma$ be a (closed) face or a union of several faces of
$G$. Define
$$H_{\Gamma}(\c;~G) = \{\v\in H(\c;~G):~\v\times\n =0 \mbox{ on }
\Gamma\}$$
and
$$ H_{\Gamma}^1(G)=\{v\in H^1(G):~v=0 \mbox{ on }
\Gamma\}.$$

\subsection{Edge and nodal element spaces}

For a polyhedron $G$, let $G$ be divided into a union of small tetrahedral
elements of size $h$, and let ${\cal T}_h$ denote the resulting triangulation of the domain
$G$. For convenience, we assume that the triangulation ${\cal T}_h$ is quasi-uniform. If the assumption is not satisfied, the estimates
built in this paper may depend on the shape parameters of elements in ${\cal T}_h$.  We use ${\cal E}_h$ and ${\cal
N}_h$ to denote the set of edges of ${\cal T}_h$ and the set of
nodes in ${\cal T}_h$ respectively. Then the Nedelec edge element
space, of the lowest order, is a subspace of piecewise linear
polynomials defined on ${\cal T}_h$: 
$$ \l{eq:nedelec}
V_h(G)=\Big\{\v\in {H}({\bf curl};~G); ~\v|_K\in R(K), ~\forall K\in\T
 \Big\},
$$
where $R(K)$ is a subset of all linear polynomials on the element
$K$ of the form:
$$ R(K)=\Big\{\a+\b\ti\x;~ \a,\b\in {\bf R}^3, ~\x\in K\Big\} \,. $$

It is known that, for any $\v\in V_h(G)$, its tangential components
are continuous on all edges in ${\cal E}_h$, and
$\v$ is uniquely determined by its moments on edges $e$ of ${\cal T}_h$:
$$ M_h(\v)=\Big\{\la_e(\v)=\int_e\v\cdot\t_e ds; ~ e\in {\cal E}_h\Big\}, $$
where $\t_e$ denotes the unit vector on an edge $e$, and this notation will be used to denote
the unit vector on any edge or union of edges, either from
an element $K\in\T$ or from $G$ itself. For example, for a face $\ff$ of $G$, the notation $\t_{\partial\ff}$ denotes
the unit vector along $\partial\ff$.

For a vector-valued function $\v$ with appropriate smoothness, we introduce its edge
element ``projection" ${\bf r}_h\v$ such that ${\bf r}_h\v\in
V_h(G)$, and ${\bf r}_h\v$ and $\v$ have the same moments as in
$M_h(\v)$. Such an operator $\r_h$ can be chosen as the standard edge
element interpolation operator or the edge
element projector ${\bf R}^1_D$ introduced in
\cite{Hiptmair-Pechstein2017}.
The operator ${\bf r}_h$ will be used in
the construction of a stable decomposition for any function $\v_h\in
V_h(G)$.

As we will see, the edge element analysis also involves frequently the
nodal element space. For this purpose, we use $Z_h(G)$ to denote the standard continuous
piecewise linear finite element space associated with the
triangulation ${\cal T}_h$.

\subsection{Some special notions and notations}

Throughout this paper, we shall frequently use the notations $\stl$.
For any two non-negative quantities $x$ and $y,$ $x\stl y$
means that $x\leq Cy$ for some constant $C$ independent of mesh size
$h$. For a non-negative quantity $z$, when $1\stl z\stl 1$, we simply write $z=O(1)$.

From now on, a ``usual" polyhedron $G$ means that $G$ is a bounded
polyhedron with a diameter $O(1)$ and vanishing first Betti number, and $G$ has a fixed number (independent of $h$) of flat faces only, but it may be non-convex.
A face $\ff$ of $G$ is always understood as a ``complete" face, i.e., $\ff$ satisfies $\ff=\pi\cap\partial G$ for the plane $\pi$ containing $\ff$, which means that $\ff$ is just a closed (flat) face instead of
only a part of a face.

We will often use {\sc f}, {\sc e} and {\sc v} to denote a general face, edge and vertex of $G$,
respectively, but use $e$ to denote a general edge in ${\cal E}_h$ associated with a triangulation ${\cal T}_h$ on $G$. Throughout this paper a face {\sc f} or an edge {\sc e} is understood as a closed subset of $\partial G$.

A connected union $\Gamma$ of several faces of $G$ is called
a {\it connected ``Lipschitz" union of some faces} if there is no isolated vertex on $\Gamma$ (see Figure 1), i.e., for any face $\ff\subset\Gamma$, there exists another face
$\ff'\subset\Gamma$ such that the intersection $\ff\cap\ff'$ is an edge; otherwise, it is called a {\it connected ``non-Lipschitz" union of some faces}. In essence,
{\it connected ``Lipschitz" union of some faces} can be understood as a Lipschitz subset of $\partial G$, but a {\it connected ``non-Lipschitz" union of some faces} is not a
Lipschitz subset of $\partial G$.

\begin{center}
\includegraphics[width=6cm,height=6cm]{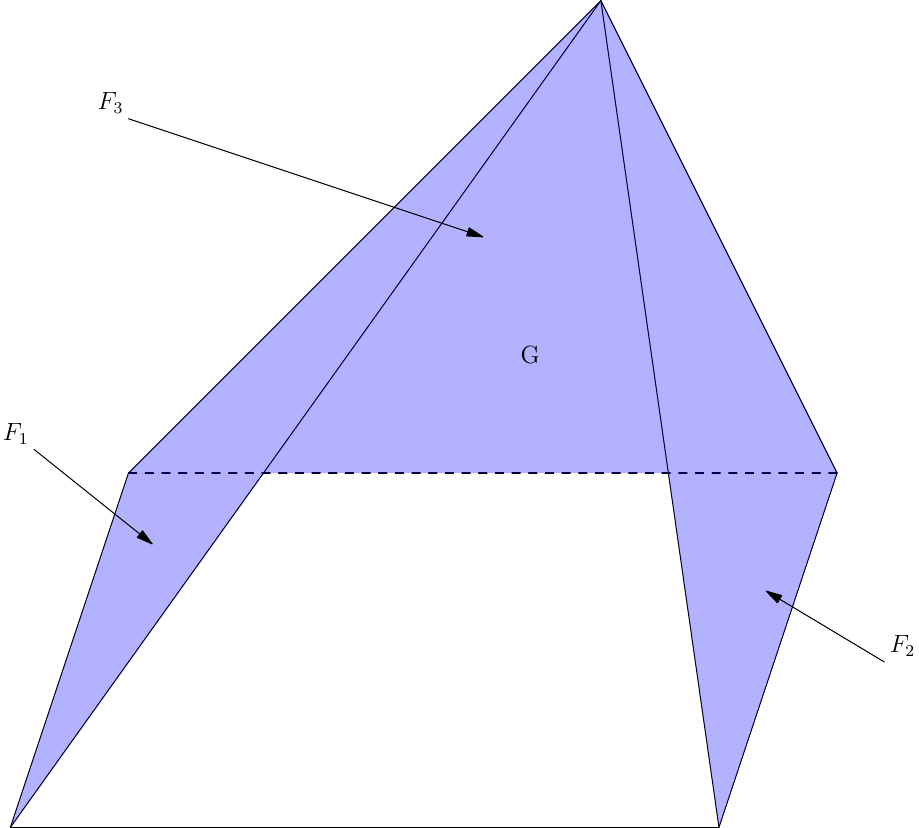}\quad\quad\includegraphics[width=6cm,height=6cm]{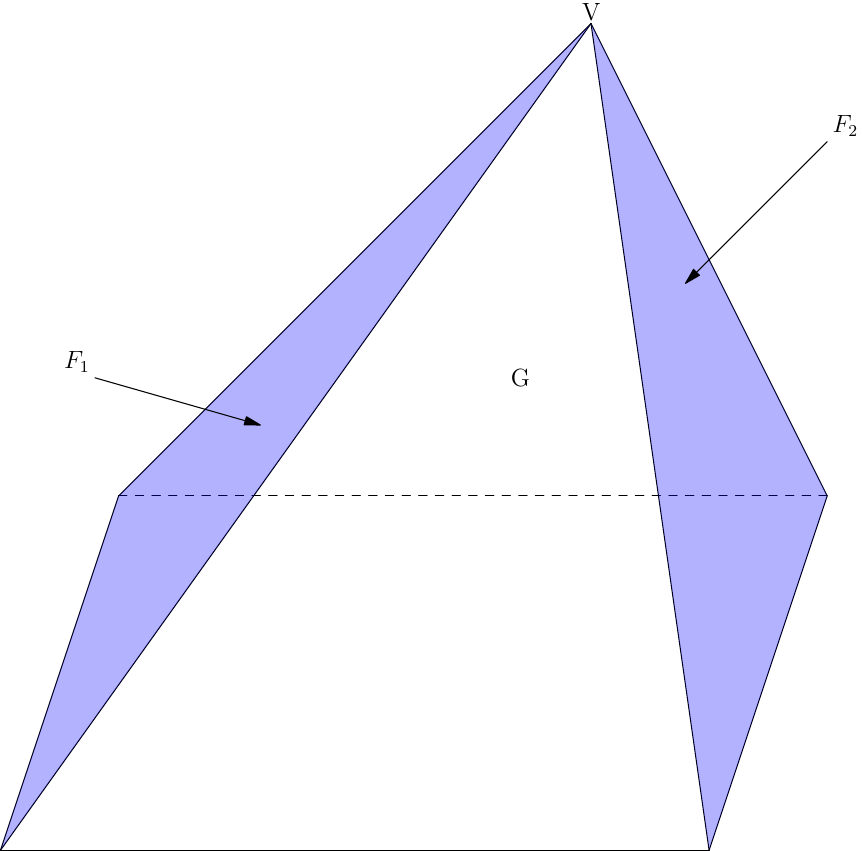}
\end{center}
Figure 1. $\ff_1$ and $\ff_2$ are two opposite lateral faces of four-sided pyramid $G$, and $\ff_3$ denotes another lateral
face that has a common edge with each of $\ff_1$ and $\ff_2$. The set $\Gamma=\ff_1\cup\ff_2\cup\ff_3$ is a {\it connected ``Lipschitz" union of three faces}~(left), but the set $\Gamma=\ff_1\cup\ff_2$
is a {\it connected ``non-Lipschitz" union of two faces}~(right), where the point $\vv$ is an isolated vertex of $\Gamma$.
\vskip 0.1in

As we will see, a connected union of several edges  has no essential difference from an edge. For simplicity of exposition, an ``edge" is always understood as an ``edge"
or a ``connected union of edges" later.

\section{Main results}

In this section we describe the new regular decompositions (and their stabilities), which will be constructed and analyzed in the upcoming sections.

Let $G$ be a usual polyhedron or a connected union of two usual polyhedra, on which a triangulation ${\mathcal T}_h$ is given, and let $\Gamma$ denote a union of faces and edges of $G$.
The following assumption will be used repeatedly:\\
{\bf Assumption 3.1}. The vector-valued function $\v_h\in V_h(G)$ has zero degrees of freedom $\la_e(\v_h)=0$ for all $e\subset\Gamma$, where $e\in {\mathcal E}_h$.

The first result can be regarded as a variant of the existing discrete regular decompositions on polyhedral domains.
\begin{theorem}~\l{teor:helm1} Let $G$ be a usual polyhedron, and let $\Gamma$ be a (may be non-connected) union of some faces and edges of $G$. Suppose that $\v_h$ satisfies {\bf Assumption 3.1}.
Then the function $\v_h$ admits a decomposition
\ee
\v_h=\nabla p_h+\r_h\w_h+{\bf R}_h \label{regulardecomposition:3.new0}
\e
for some $p_h\in Z_h(G)$, $\w_h\in (Z_h(G))^3$ and
${\bf R}_h\in V_h(G)$ such that $p_h=0$ and $\w_h=\0$ {on}
${\Gamma}$, and $\la_e({\bf R}_h)=0$ for all
$e\subset\Gamma$. Moreover, we have
\ee \|\w_h\|_{1,G}+h^{-1}\|{\bf R}_h\|_{0,G}\stl\log(1/h)\|\v_h\|_{\c,G}\label{stab:4.new19}
\e
and
\ee
\|\w_h\|_{0,G}+\|p_h\|_{1,G}\stl\log(1/h)\|\v_h\|_{\c,G}. \label{stab:4.new20}
\e
\end{theorem}

We would like to emphasize that, in Theorem \ref{teor:helm1}, the union set $\Gamma$ may contain some edges of the polyhedron, which is different from existing results.
\begin{remark} In particular, when $\Gamma$ is a connected set of $\partial G$,
the norm on the right side of (\ref{stab:4.new19}) can be replaced by the $\c$ semi-norm. If there is not an isolated edge in $\Gamma$, i.e., $\Gamma$ is a union of faces only,
then the norm on the right side of (\ref{stab:4.new20}) can be replaced by the $L^2$ norm. When $\Gamma=\cup_{j=1}^J\Gamma_j$ with $\Gamma_j$ being a connected
``Lipschitz" union of some faces of $G$, the factor $\log(1/h)$ in the estimates (\ref{stab:4.new19}) and (\ref{stab:4.new20}) can be dropped.
\end{remark}

We would like to extend the above theorem to the case of non-Lipschitz polyhedra.
Let $G_1$ and $G_2$ be two intersecting (usual) polyhedra, and define $G$ as their union: $\bar{G}=\bar{G}_1\cup \bar{G}_2$. We consider two
 particular cases: (1) the intersection $\bar{G}_1\cap \bar{G}_2$ is an edge $\E$ of $G_1$ and $G_2$; (2) the intersection $\bar{G}_1\cap \bar{G}_2$ is
a vertex $\vv$ of $G_1$ and $G_2$. For the two cases, $G$ is not a Lipschitz domain.

 We assume that the two triangulations on $G_1$ and $G_2$ are matching on the intersection $\bar{G}_1\cap \bar{G}_2$, i.e.,
 their nodes are coincident on $\bar{G}_1\cap \bar{G}_2$. An edge finite element function $\v_h\in V_h(G)$
 means that $\v_h|_{G_i}\in V_h(G_i)$ ($i=1,2$) and $\v_h$ has continuous tangential complements on $\E$ for Case (1). A nodal finite element function $p_h\in Z_h(G)$
 means that $p_h|_{G_i}\in Z_h(G_i)$ ($i=1,2$) and $p_h$ is continuous on the intersection $\bar{G}_1\cap \bar{G}_2$. Norms of finite element functions on $G$ are understood
 in the natural manner, for instance, the norm of $\v_h\in V_h(G)$ is defined as $\|\v_h\|_{\c, G}=(\|\v_h\|^2_{\c, G_1}+\|\v_h\|^2_{\c, G_2})^{{1\over 2}}$. For a union $\Gamma$
 of faces and edges of $G$, define $\Gamma_i=\Gamma\cap\partial G_i$ ($i=1,2$).


The following result is a direct extension of Theorem \ref{teor:helm1} to Case (1).

\begin{theorem}\l{thm:thm5.1} Let $G$ be a union of two usual polyhedra $G_1$ and $G_2$, with $\bar{G}_1\cap \bar{G}_2$ being the common edge $\E$ of $G_1$ and $G_2$,
and let $\Gamma$ be a union of some faces and edges of $G$. Suppose that $\v_h$ satisfies {\bf Assumption 3.1}. Then $\v_h$ admits a decomposition
\ee \v_h=\nabla p_h+\r_h\w_h+{\bf R}_h \label{5.new1}
\e
for some $p_h\in Z_h(G)$, $\w_h\in (Z_h(G))^3$ and
${\bf R}_h\in V_h(G)$ such that $p_h=0$ and $\w_h=\0$ {on}
${\Gamma}$, and $\la_e({\bf R}_h)=\0$ for all $e\subset\Gamma$. Moreover, we have
\ee
\|\w_h\|_{1,G}+h^{-1}\|{\bf R}_h\|_{0,G}\stl\log(1/h)\|\v_h\|_{\c,G}\label{stab:5.new1}
\e
and
\ee
\|\w_h\|_{0,G}+\|p_h\|_{1,G}\stl\log(1/h)\|\v_h\|_{\c,G}. \label{5.new2}
\e
In particular, when $\Gamma$ contains faces only and $\E\subset \Gamma_i$ for $i=1,2$, the logarithmic factor in the above estimates can be dropped. If
$\Gamma$ is connected, then the complete norm on the right side of (\ref{stab:5.new1}) can be replaced by the $\c$ semi-norm.
\end{theorem}

Now we consider Case (2), i.e., $\bar{G}_1\cap\bar{G}_2=\vv$ is a vertex of $G_1$ and $G_2$. For this case, the problem is a bit complicated, and
a similar result with Theorem \ref{thm:thm5.1} cannot be built yet except that the considered function $\v_h$ meets an additional condition.
\begin{theorem}\l{thm:thm5.2} Let $G$ be a union of two usual polyhedra $G_1$ and $G_2$, with $\bar{G}_1\cap \bar{G}_2$ being the common vertex $\vv$ of $G_1$ and $G_2$.
Let $\Gamma$ denote a union of some faces and edges of $G_1$ and $G_2$. Suppose that $\v_h$ satisfies {\bf Assumption 3.1}.
Then there exists a functional ${\cal F}$ (whose definition will be given in Section 6) such that, if $\v_h$ satisfies the constraint
${\cal F}\v_h=0$, the function $\v_h$ has a decomposition
\ee \v_h=\nabla p_h+\r_h\w_h+{\bf R}_h \label{5.new8}
\e
for some $p_h\in Z_h(G)$, $\w_h\in (Z_h(G))^3$ and
${\bf R}_h\in V_h(G)$ such that $p_h$, $\w_h$ and ${\bf R}_h$ have zero degrees of freedom on
${\Gamma}$. Moreover, we have
\ee
\|\w_h\|_{1,G}+h^{-1}\|{\bf R}_h\|_{0,G}\stl\log(1/h)\|\v_h\|_{\c,G}
\label{stab:5.new4}\e
and
\ee
\|\w_h\|_{0,G}+\|p_h\|_{1,G}\stl\log(1/h)\|\v_h\|_{\c,G}. \label{5.new9}
\e
When both $\Gamma_1$ and $\Gamma_2$ are connected, the norm on the right side of (\ref{stab:5.new4})
can be replaced by the $\c$ semi-norm.
\end{theorem}
\begin{remark} The additional condition ${\cal F}\v_h=0$ in Theorem \ref{thm:thm5.2} seems absolutely necessary.
The value of an edge finite element function $\v_h$ is not uniquely defined at the common vertex $\vv$ of $G_1$ and $G_2$,
but the two nodal finite element functions defined by the regular decomposition
of $\v_h$ are required to be continuous at the vertex $\vv$. Intuitively understanding, there exists a gap between $\v_h$
and the gradient of the scalar finite element function $p_h$. The role of the additional condition on $\v_h$ is to fill this gap.
\end{remark}

We can also consider the case that $G$ is a union of more polyhedrons. Let $G_1,G_2,\cdots,G_s$ be usual polyhedrons that may be non-convex, and let $G$ be
a union of $G_1,G_2,\cdots,G_s$ (then $G$ is a non-Lipchitz domain) such that the intersection of any two polyhedrons in $G_1,G_2,\cdots,G_s$ just is the same vertex $\vv$,
i.e., $\bar{G}_i\cap\bar{G}_j=\vv$ for $i\not=j$, which implies that $\bar{G}_1\cap\bar{G}_2\cap\cdots\cap\bar{G}_s=\vv$.

\begin{theorem}\l{thm:thm5.4} Let $G_i$ ($i=1,\cdots,s;~s\geq 3$) and $G$ be defined as above, and let $\Gamma$ denote a union of some faces and edges of $G_1,\cdots,G_s$.
Suppose that $\v_h$ satisfies {\bf Assumption 3.1}. Then there exist $s-1$ functionals ${\cal F}_i$  (whose definition will be given at the ending of Section 6) such that, if $\v_h$ satisfies the constraints
${\cal F}_i\v_h=0$ for $i=1,\cdots,s-1$, the function $\v_h$ admits
a regular decomposition like (\ref{5.new8}) and the resulting functions satisfy the conditions and estimates described in Theorem 3.3.
\end{theorem}

\section{Basic regular decompositions}
\l{sec:variant2} \setcounter{equation}{0}
In this section we develop regular decompositions for vector-valued functions that have zero tangential trace on a part of the boundary of a Lipschitz polyhedron.

A regular decomposition is called {\it $\c$-bounded} if the $H^1$ norms of the two $H^1$ functions defined by the decomposition can be controlled
by the $\c$ semi-norm of the original vector-valued function;
a regular decomposition is called {\it $L^2$-stable} if the $L^2$ norm of the vector-valued $H^1$ function and the $H^1$ norm of the scalar $H^1$ function can be controlled
by the $L^2$ norm of the original vector-valued function.

In this section we first build $\c$-bounded regular decompositions on a class of Lipschitz polyhedrons and then develop $\c$-bounded and $L^2$-stable
regular decompositions on a class of non-convex polyhedrons .
\subsection{${\bf Curl}$-bounded regular decompositions on Lipschitz polyhedrons}

Let $G$ be a bounded Lipschitz polyhedron with a diameter $O(1)$ and vanishing first Betti number.
Assume that $G$ is a union of tetrahedral elements, associated with a triangulation ${\mathcal T}_h$ ($h\ll 1$), with $h$ being the mesh width.
For each $r$ ($1\leq r\leq J$), let $\Gamma_r\subset\partial G$ be a connected union
of faces of some elements of ${\mathcal T}_h$ (see Figure 2).\\
{\bf Assumption 4.1}. For each $\Gamma_r$, there exists a bounded Lipschitz domain $D_r$ with a diameter $O(1)$ and vanishing first Betti number,
such that (i) $\bar{D}_r\cap\bar{G}=\Gamma_r$; (ii) $dist(D_i,D_j)\geq\delta$ for any $i\not=j$, where $\delta$ is a positive constant;
(iii) $B=G\bigcup(\cup_{r=1}^J\bar{D}_r)$ is a Lipschitz domain with vanishing first Betti number.

The condition (ii) implies that $dist(\Gamma_i,\Gamma_j)\geq\delta$ for any $i\not=j$. This assumption can be viewed as an extension of the assumption given in \cite{Hiptmair-Zeng2009} and
is essentially equivalent to the one introduced in \cite{Hiptmair-Pechstein2017}, but
it seems more direct. Before presenting the main result of this subsection, we would like to give two examples satisfying Assumption 4.1.

\begin{center}
\includegraphics[width=6cm,height=6cm]{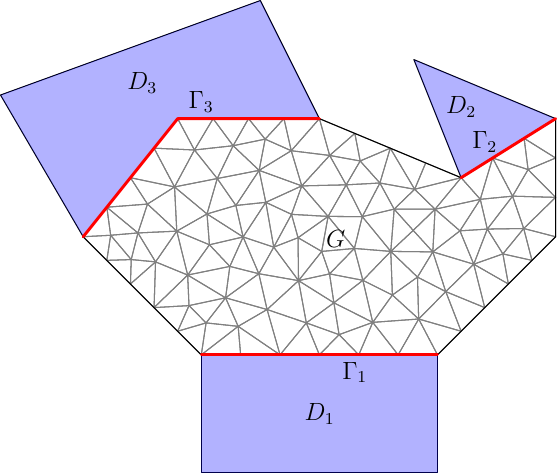}
\end{center}
\vskip 0.1in
Figure 2: The cross-section of the domains described in Assumption 4.1, where the red parts denote the connected components of $\Gamma$, the blue parts denote
the auxiliary subdomains.
\vskip 0.1in

{\bf Example 4.1}. Let domain $G$ be a usual polyhedron with fixed number of (plane) faces only, and let each $\Gamma_r$ be a complete face or a connected union of several complete faces. Assume that $\bar{\Gamma}_i\cap\bar{\Gamma}_j=\emptyset$ for $i\not=j$.

We will explain that Assumption 4.1 is met as long as each $\Gamma_r$ is a {\it connected ``Lipschitz" union of some faces}, which means that $\partial\Gamma_r$ is connected and closed but is not self intersecting.
We only need to define a suitable domain $D_r$ for each $\Gamma_r$. It is easy to see that $\partial\Gamma_r$ can be written as $\partial\Gamma_r=\cup_{l=1}^{n_r}\E^{(r)}_l$, where $\E^{(r)}_l$ is an edge of
some face $\ff^{(r)}_l\subset\Gamma_r$ ($l=1,\cdots,n_r$). For each $\E^{(r)}_l$, we choose a half plane $\pi^{(r)}_l$ such that $\pi^{(r)}_l\cap\bar{G}=\E^{(r)}_l$ and the interfacial
angle between $\pi^{(r)}_l$ and the face $\ff^{(r)}_l$ is an acute angle (this means that the boundary line of $\pi^{(r)}_l$ contains $\E^{(r)}_l$ and $\pi^{(r)}_l$ extends to the exterior of $G$).
For a convex polyhedron $G$, $\pi^{(r)}_l$ can be simply chosen as the extending half plane of the face that has the common edge $\E_l^{(r)}$ with $\ff_l^{(r)}$ (see the left graph in Figure 3). If $\Gamma_r$ and all the half planes $\pi^{(r)}_1,\cdots,\pi^{(r)}_{n_r}$ just
surround a polyhedron on one side of $G$, then this polyhedron can be chosen as $D_r$ and it satisfies $\bar{D}_r\cap\bar{G}=\Gamma_r$. Otherwise, we can add one (or several) suitable plane $\pi^*$ such that
$\Gamma_r$, $\pi^*$ and all the half planes $\pi^{(r)}_1,\cdots,\pi^{(r)}_{n_r}$ surround a polyhedron on one side of $G$, and we choose this polyhedron as $D_r$.
It is easy to see that Assumption 4.1 is met for this choice of $D_r$. In particular, when $G$ is convex or $\Gamma$ contains all the concave parts of $\partial G$ (its precise definition will be given later), the resulting extended domain $B$ is convex.

However, if some $\Gamma_r$ is a connected ``non-Lipschitz" union of some faces of $G$ (for example, $G$ is
a four-sided pyramid and $\Gamma_r$ is a union of two opposite lateral faces of $G$, see the right graph of Figure 1), then
$\partial\Gamma_r$ is self intersecting and so there does not exist a Lipschitz domain $D_r$ satisfying Assumption 4.1.

\begin{center}
\includegraphics[width=6cm,height=6cm]{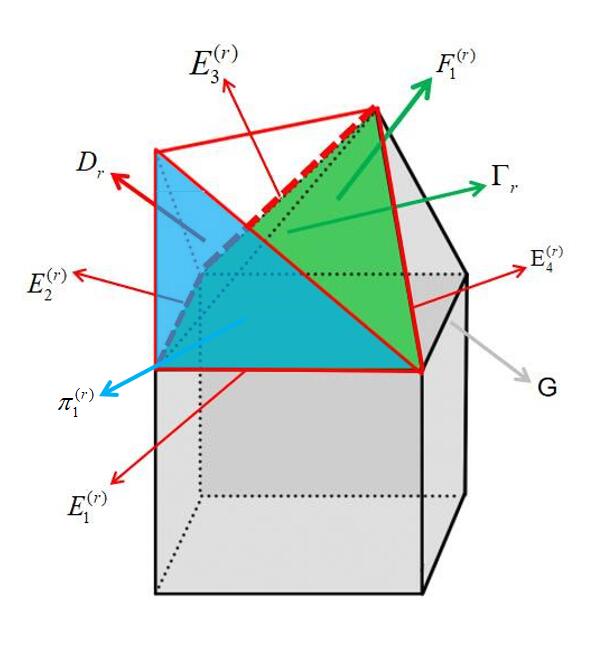} \quad \quad\quad \includegraphics[width=5cm,height=5cm]{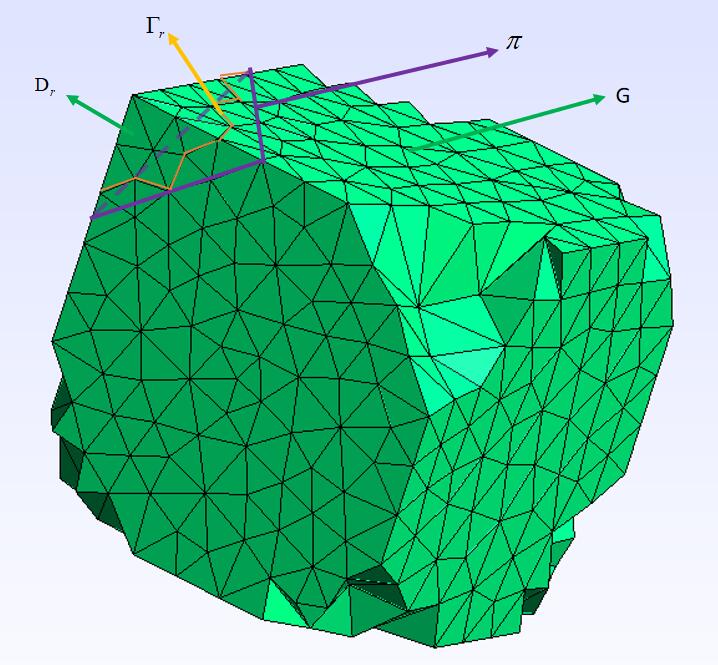}
\end{center}
\vskip 0.1in
Figure 3: In the left graph, $G$ denote a union of the cube and the four-sided pyramid, $\Gamma_r$ denotes a union of two neighboring triangle
faces of the four-sided pyramid, the polyhedron $D_r$ is generated by four half planes that extend to the exterior of $G$ from the red edges $\E^{(r)}_1,\cdots,\E^{(r)}_4$
(the blue triangle face of $D_r$ is a part of the half plane $\pi^{(r)}_1$ that contains the edge $\E^{(r)}_1$ and extends above).
In the right graph, $G$ is generated from a cube domain $\Omega$ in the manner described in Example 4.2, $D_r$ is the Lipschitz polyhedron at the upper-left corner
and with yellow edges (the tetrahedron with purple edges denotes the cut domain $\hat{D}_r$, $\pi$ denote the cut plane), the global green polyhedron is a union of $G$ and one $D_r$.
\vskip 0.1in

{\bf Example 4.2}: For a usual polyhedron $\Omega$, let $\Omega^{\partial}\subset\Omega$ be the polyhedron generated by cutting its convex angles with some planes,
such that every cut domain $\hat{D}_r$ has much larger size than $h$ and has a suitable distance from the other cut domains. Notice that the polyhedron $\Omega^{\partial}$ may be not a union of elements.
Let $G$ be a union of all the elements, each of which has at least one vertex in $\Omega^{\partial}$ (see the right graph in Figure 3). Then $G$ is a Lipschitz polyhedron, which can be regarded as a small perturbation of the
usual polyhedron $\Omega^{\partial}$. Let $D_r$ be a union of all the elements contained in $\hat{D}_r$ and set $\Gamma_r=\bar{G}\cap\bar{D}_r$.
We define $B=G\bigcup(\cup_{l=1}^J\bar{D}_r)=\Omega$. It is clear that Assumption 4.1 is met for this example.

Set $\Gamma=\cup_{l=1}^J\Gamma_l$. We use the same notations introduced in Section 2.
The following result on regular decompositions can be found in \cite{Hiptmair-Pechstein2017}.

\begin{lemma}\label{regular decomposition} Let Assumption 4.1 be satisfied.
    For any $\v\in H_{\Gamma}(\c;~G)$, there exists a vector-valued function $\Phi\in
    (H_{\Gamma}^1(G))^3$ and a scalar function $p\in H_{\Gamma}^1(G)$ such that
\ee \v=\Phi+\nabla p. \l{regul}\e
Moreover, we have \ee
|\Phi|_{1,G}\stl\|\v\|_{\c,G}~~\mbox{and}~~
\|\Phi\|_{0,G}+\|p\|_{1,G}\stl\|\v\|_{0,G}. \label{stab4.0}\e In
particular, if $\v=\v_h\in V_h(G)\cap H_{\Gamma}(\c;~G)$, then there exist $p_h\in
Z_h(G)\cap H_{\Gamma}^1(G)$, $\Phi_h\in (Z_h(G)\cap
H^1_{\Gamma}(G))^3$ and ${\bf R}_h\in V_h(G)\cap H_{\Gamma}(\c;~G)$ such that
\ee \v_h=\r_h\Phi_h+\nabla p_h+{\bf R}_h, \label{decom4.00}\e
where $\r_h$ denotes the edge element projector ${\bf R}^1_D$ introduced in \cite{Hiptmair-Pechstein2017}. Moreover, they satisfy the stabilities
\ee |\Phi_h|_{1,G}+h^{-1}\|{\bf
R}_h\|_{0,G}\stl\|\v\|_{\c,G}~~~\mbox{and}~~\|\Phi\|_{0,G}+\|p_h\|_{1,G}\stl\|\v\|_{0,G}.\label{stab4.02}\e
\end{lemma}
\hfill $\Box$

In the following lemma, we present {\it $\c$-bounded} regular decompositions, which have different stability estimates from that in Lemma \ref{regular decomposition}.
\begin{lemma}\label{regular decomposition-3.new} Let $J=1$, i.e., $\Gamma$ be a connected union of closed faces of elements, and let Assumption 4.1 be satisfied.
For any $\v\in H_{\Gamma}(\c;~G)$, there exists a vector-valued function $\Phi\in(H_{\Gamma}^1(G))^3$ and a scalar function $p\in H_{\Gamma}^1(G)$ satisfying (\ref{regul}),
such that
\ee
\|\Phi\|_{1,G}\stl\|\c\,\v\|_{0,G}~~\mbox{and}~~
\|p\|_{1,G}\stl\|\v\|_{\c,G}.\label{stab4.0new}
\e
When $\v=\v_h\in V_h(G)\cap H_{\Gamma}(\c;~G)$, there exist $p_h\in
Z_h(G)\cap H_{\Gamma}^1(G)$, $\Phi_h\in (Z_h(G)\cap
H^1_{\Gamma}(G))^3$ and ${\bf R}_h\in V_h(G)\cap H_{\Gamma}(\c;~G)$,
such that they satisfy the decomposition (\ref{decom4.00}) with $\r_h$ being the standard edge element interpolation, and they have the stability estimates
\ee\|\Phi_h\|_{1,G}+h^{-1}\|{\bf R}_h\|_{0,G}\stl\|\c\,\v_h\|_{0,G}~~
\mbox{and}~~\|p_h\|_{1,G}\stl\|\v_h\|_{\c,G}.
\label{stab4.01new}\e
In particular, if the extended domain $B$ is convex, the second inequality in (\ref{stab4.0new}) (and (\ref{stab4.01new})) can be replaced by
\ee
\|\Phi\|_{0,G}+\|p\|_{1,G}\stl\|\v\|_{0,G}\quad\mbox{and}\quad
\|\Phi_h\|_{0,G}+\|p_h\|_{1,G}\stl\|\v_h\|_{0,G}.\label{stab-new00}
\e
\end{lemma}
{\it Proof}. The basic ideas follow the arguments in \cite{Regul,PasciakZhao2002}, here the proof seems more direct. When $J=1$, we have $\Gamma=\Gamma_1$ and $D=D_1$.
Let $B$ be the extended domain defined in Assumption 4.1.
We extend $\v$ onto the global $B$ by zero, i.e., the extension $\tilde{\v}$ satisfying
$\tilde{\v}=\v$ on $G$ and $\tilde{\v}=\0$ on $\bar{D}$.
Since $\v\times{\bf n}=\0$ on $\Gamma$, we have $\tilde{\v}\in
H(\c;~B)$. By the regular decomposition (Lemma 2.4 in \cite{hipt2002}), we get\footnote{Let ${\bf N}_{{\mathcal H}}$ and ${\bf R}$ be the mappings defined in Lemma 2.4 of \cite{hipt2002}.
We define $\w={\bf R}\tilde{\v}$. Since the first Betti number of $B$ vanishes and so we have ${\bf N}_{{\mathcal H}}\tilde{\v}=\0$.
The property ${\bf R}|_{{\mathcal {\bf H}}(\bf{\c} \0; B)}=\0$ means that
$\|\w\|_{1,B}\stl \|\c~\tilde{\v}\|_{0,B}$. }
\ee
\tilde{\v}=\w+\nabla \varphi~~~~\mbox{on}~B,\label{Helm5.new1} \e
with $\w\in (H^1(B))^3$ and $\varphi\in H^1(B)$. Moreover, $\w$ and $\varphi$ satisfy
\ee
\|\w\|_{1,B}\stl\|\c\,\tilde{\v}\|_{0,B}=\|\c\,\v\|_{0,
G}~~\mbox{and}~~
\|\nabla\varphi\|_{0,B}\stl\|\tilde{\v}\|_{\c,B}=\|\v\|_{\c,G}.
\label{stab5.new1}\e
Here we require that $\varphi$ has zero average value on $D$ (not on $B$), which is written as $\gamma_D(\varphi)=0$.
When $B$ is convex, by the orthogonal Helmholtz decomposition in \cite{Gir1986}, we have \ee
\|\w\|_{1,B}\stl\|\c\,\v\|_{0, G}~~\mbox{and}~~
\|\w\|_{0,B}+\|\nabla\varphi\|_{0,B}\stl\|\v\|_{0,G}.
\label{stab5.new01}\e
Notice that $\tilde{\v}=\0$ on $\bar{D}$, we
have $\nabla\varphi=-\w$ on $D$, and so $\varphi\in H^2(D)$. Let
$\tilde{\varphi}\in H^2(B)$ be a stable extension of $\varphi$
from $D$ onto the global $B$. It follows by (\ref{Helm5.new1}) that
\ee \tilde{\v}=(\w+\nabla\tilde{\varphi})+\nabla
(\varphi-\tilde{\varphi})~~~~\mbox{on}~B.\label{Helm5.new2} \e
Define $\Phi=\w+\nabla\tilde{\varphi}$ and
$p=\varphi-\tilde{\varphi}$. Then we have $\Phi\in (H^1(B))^3$ and $p\in
H^1(B)$, and they satisfy (\ref{regul}). Since $\tilde{\varphi}=\varphi$ on $\bar{D}$, we have
$p=0$ on $\bar{D}$, which means that $p\in H_{\Gamma}^1(G)$ and $\nabla p=\0$ on $\bar{D}$. Therefore
$\Phi=\0$ on $\bar{D}$ and so $\Phi=\0$ on $\Gamma$, which implies that
$\Phi\in (H_{\Gamma}^1(G))^3$.

By the definition of $\Phi$ and the first inequality in
(\ref{stab5.new1}), together with Friedrichs inequality and the relation $\nabla\varphi=-\w$ on $D$, we get (notice that $D$ is simply-connected,
so Friedrichs inequality can be used on $D$)
\beq \|\Phi\|_{1,G}&\leq&
\|\w\|_{1,G}+\|\tilde{\varphi}\|_{2,G}\stl
\|\c\,\v\|_{0,G}+\|\varphi\|_{2,D}\cr
&\stl&\|\c\,\v\|_{0,G}+\|\nabla\varphi\|_{1,D}
\stl\|\c\,\v\|_{0,G}+\|\w\|_{1,D}
\cr&\stl&\|\c\,\v\|_{0,G}+\|\w\|_{1,B}\stl\|\c\,\v\|_{0,G}.\label{inequality-3.new00}
\eq
We can further obtain the second inequality of (\ref{stab4.0new})
by (\ref{regul}) since $p=0$ on $\Gamma$.

Let ${\bf Q}_{h}: (H^1_{\Gamma}(G))^3\rightarrow (Z_h(G)\cap
H^1_{\Gamma}(G))^3$ denote the Scott-Zhang interpolation (see \cite{ScottZ1990}), which can preserve Dirichlet zero boundary condition, and define
$\Phi_h={\bf Q}_{h}\Phi$. In addition,  let $\r_h$ be the standard edge element interpolation operator. By Lemma 4.7 in \cite{Amr1998} (notice that $\c~\Phi=\c~\v_h$), we know that the interpolation functions $\r_h\Phi$ and $\r_h\Phi_h$ are well defined.
Define
$${\bf R}_h=\r_h(\Phi-\Phi_h)=({\bf I}-{\bf Q}_{h})\Phi+(\r_h-{\bf I})(\Phi-\Phi_h).$$
Then
\ee
\v_h=\r_h\Phi+\r_h(\nabla p)=\r_h\Phi+\nabla p_h=\r_h\Phi_h+\nabla p_h+{\bf R}_h\label{decomp-3.new00}
\e
with $p_h\in (Z_h(G)\cap H^1_{\Gamma}(G))^3$.  By the $H^1$ stability of the operator ${\bf Q}_h$ and the first estimate in (\ref{stab4.0new}), we directly
get
\ee
\|\Phi_h\|_{1,G}\stl\|\Phi\|_{1,G}\stl\|\c~\v_h\|_{0,G}.\label{stab4.01new1}
\e
From the proof of Lemma 4.3 in \cite{HZ2}, we can see that ($\c~(\Phi-\Phi_h)=\c~\v_h-\c~\Phi_h$ is a piecewise constant vector)
$$ \|(\r_h-{\bf I})(\Phi-\Phi_h)\|_{0,G}\stl h(|\Phi-\Phi_h|_{1,G}+\|\c~(\Phi-\Phi_h)\|_{0,G})\stl h(|\Phi|_{1,G}+|\Phi_h|_{1,G}). $$
This, together with the first inequalities in (\ref{stab4.0new})-(\ref{stab4.01new}), yields
$$ \|(\r_h-{\bf I})(\Phi-\Phi_h)\|_{0,G}\stl h\|\c~\v_h\|_{0,G}. $$
Then, by the $L^2$ approximation estimates for ${\bf Q}_h$ (see Theorem 4.1 of \cite{ScottZ1990}), we further obtain
\beqx
\|{\bf R}_h\|_{0,G}&\leq&\|({\bf I}-{\bf Q}_{h})\Phi\|_{0,G}+\|(\r_h-{\bf I})(\Phi-\Phi_h)\|_{0,G}\cr
&\stl&h|\Phi|_{1,G}+h\|\c~\v_h\|_{0,G}\stl h\|\c~\v_h\|_{0,G}.
\eqx
This, together with (\ref{stab4.01new1}), gives the first inequality in (\ref{stab4.01new}). The second inequality in (\ref{stab4.01new}) can be derived immediately by (\ref{decomp-3.new00}).

If $B$ is convex, then we can use the second inequality in (\ref{stab5.new01}) to verify the first estimate in (\ref{stab-new00}) in an analogous way with (\ref{inequality-3.new00}).
Using the $L^2$ approximation of ${\bf Q}_h$ and the first estimate in (\ref{stab-new00}) yields
\ee
 \|\Phi_h\|_{0,G}\stl\|\Phi\|_{0,G}+h|\Phi|_{1,G}\stl\|\Phi\|_{0,G}+h|\c~\v_h|_{1,G}\stl\|\Phi\|_{0,G}+\|\v_h\|_{0,G}\stl\|\v_h\|_{0,G}.
\label{stab-lemma3.2-new}\e
Similarly, using the equality (\ref{decomp-3.new00}) and the $L^2$ approximation estimates for ${\bf r}_h$, together with the first estimate in (\ref{stab-new00}), leads to (notice that
$\c~\Phi=\c~\v_h$)
\beqx
\|\nabla p_h\|_{0,G}&\leq&\|\v_h\|_{0,G}+\|\r_h\Phi\|_{0,G}\stl\|\v_h\|_{0,G}+\|\Phi\|_{0,G}+\|(\r_h-{\bf I})\Phi\|_{0,G}\cr
&\stl&\|\v_h\|_{0,G}+h(|\Phi|_{1,G}+\|\c~\Phi\|_{0,G})\cr&\stl&\|\v_h\|_{0,G}+h\|\c~\v_h\|_{0,G}\stl\|\v_h\|_{0,G}.
\eqx
This, together with (\ref{stab-lemma3.2-new}), gives the second inequality in (\ref{stab-new00}).

We would like to emphasize that all the constants omitted in the proof are independent of $h$,
since the involved domains have Lipschitz constants independent
of $h$. \hfill $\Box$

\begin{remark} It can be seen from Lemma 4.1 and Lemma 4.2 that the functions defined by different discrete regular decompositions may have different
stabilities (since the functions are constructed in different manners). In general, the {\it $\c$-bounded} stability
\ee
\|\Phi_h\|_{1,G}+h^{-1}\|{\bf R}_h\|_{0,G}\stl\|\c\,\v_h\|_{0,G}\label{estimat-3.new001}
\e
and the $L^2$ stability
\ee \|\Phi_h\|_{0,G}+\|p_h\|_{1,G}\stl\|\v_h\|_{0,G}\label{stab-3.new001}
\e
cannot be satisfied simultaneously. The inequality (\ref{estimat-3.new001}) implies that a $\c$-free vector-valued function must be
the gradient of some scalar $H^1$ function.
Lemma 4.2 gives some conditions that can guarantee both the estimates to be held. It is clear that, when $\Gamma=\partial G$ (which is connected), the conditions are met,
where the extended domain $B$ can be chosen as a ball (for this case, the result has essentially been built in \cite{PasciakZhao2002}). Other situations satisfying the conditions
was described in Example 4.1.
As we will see in \cite{Hu2-2017}, when higher order term in the considered Maxwell equations is dominant,
the semi-norm controlled estimate (\ref{estimat-3.new001}) will play a key role in the analysis. Unfortunately, if $\Gamma$ is non-connected (i.e., $J\geq 2$),
this estimate is false (a counterexample can be constructed). In fact, the domain $D$ is non-connected for this case, so Friedrichs inequality
cannot be used to prove (\ref{inequality-3.new00}). By the way, if we require that the function $\varphi$ satisfying (\ref{Helm5.new1}) has zero average value on $B$,
which is a natural idea, the semi-norm controlled estimate (\ref{inequality-3.new00}) cannot be obtained.
\end{remark}

\subsection{${\bf Curl}$-bounded and $L^2$-stable regular decompositions on non-convex polyhedrons}

In this subsection we establish slightly weaker results than (\ref{estimat-3.new001})-(\ref{stab-3.new001}) for a class of non-convex polyhedrons.
To this end, we first give several basic auxiliary results. \\
{\bf Proposition 4.1}. Let $\Omega$ be a usual polyhedron, and assume that $\w_h\in
(Z_h(\Omega))^3$. Then we have $\c~(\r_h\w_h)=\c~\w_h$ and $\|\r_h\w_h\|_{0,\Omega}
\stl\|\w_h\|_{0,\Omega}$.\\
\no{\it Proof}. Let $W_h(\Omega)$ denote the Raviart-Thomas finite
element space of the lowest order, and let $\Pi_h$ be the
interpolation operator into $W_h(\Omega)$. Since $\w_h\in (Z_h(\Omega))^3$, we
have $\c~\w_{h}\in W_h(\omega)$. Then
$$\c~(\r_h\w_{h})=\Pi_h\c~\w_{h}=\c~\w_{h}.
$$
The desired inequality can be derived by the approximation property of $\r_h$
and the inverse estimate of finite element functions.
\hfill
$\Box$

The following two results come from the so called ``edge" lemma and ``face" lemma (cf. Lemma 4.16 and Lemma 4.24 in \cite{Wid05}, or Lemma 4.9-Lemma 4.10 in \cite{s13}). Here the considered domain $\Omega$ has a 
diameter $O(1)$ and the space $H^{{1\over 2}}_{00}(\ff)$ is defined in the standard manner (cf. \cite{s13}).
\begin{lemma} Let $\Omega$ be a usual polyhedron, and $\E$ be an edge of $\Omega$. Then, for $v_h\in Z_h(\Omega)$, we have
\ee
\|v_h\|_{0,\E}\stl \log^{{1\over 2}}(1/h)\|v_h\|_{1, \Omega}. \label{edge-lemma}
\e
\end{lemma}
\begin{lemma} Let $\Omega$ be a usual polyhedron, and $\ff$ be a face of $\Omega$. For $v_h\in Z_h(\Omega)$, define $I^0_{\ff}v_h\in V_h(\partial\Omega)$ such that $I^0_{\ff}v_h$ equals
$v_h$ at the nodes in the interior of $\ff$ and vanishes at all the other nodes on $\partial\Omega$. Then
\ee
\|I^0_{\ff}v_h\|_{H^{{1\over 2}}_{00}(\ff)}\stl \log(1/h)\|v_h\|_{1, \Omega}. \label{face-lemma}
\e
\end{lemma}

By Lemma 4.24 in \cite{Wid05}, we can easily derive the following result
\begin{lemma} Let $\Omega$ be a usual polyhedron, and $\ff$ be a face of $\Omega$. Assume that $v_h\in Z_h(\Omega)$. Then there exists a function $v_h^{\ff}\in Z_h(\Omega)$
such that $v_h^{\ff}$ equals $v_h$ at all the nodes in the interior of $\ff$
and vanishes on $\partial \Omega\backslash\ff$. Moreover, the function $v_h^{\ff}$ satisfies the estimates
\ee
\|v^{\ff}_h\|_{1, \Omega}\stl \log(1/h)\|v_h\|_{1, \Omega}. \label{face-extension1}
\e
and
\ee
\|v^{\ff}_h\|_{0, \Omega}\stl \|v_h\|_{0, \Omega}. \label{face-extension2}
\e
\end{lemma}

The inequality (\ref{face-extension1}) is included in the results of Lemma 4.24 in \cite{Wid05}, and
the inequality (\ref{face-extension2}) can be directly obtained by the definition of the function $v^{\ff}_h$, which depends on the weight function $\vartheta_{\ff}$ constructed in Lemma 3.3.6 of \cite{Cas1996} and Lemma 4.23 of \cite{Wid05}.

Let $G$ be a usual non-convex polyhedron, which is a non-overlapping union of
several convex polyhedra $G_1,\cdots,G_{m_0}$. As usual, we assume that each $G_r$ is a union of elements and the intersection of any two neighboring polyhedra $G_i$ and $G_j$ is either their common face, or their
common edge, or their common vertex, i.e., the intersection cannot be only a part of one face or one edge of $G_i$ or $G_j$.

We need to give precise definitions of convex ``complete" face and ``concave part" of $\partial G$.
Let $\ff\subset \partial G$ be a connected union of some triangles, and let $\pi$ be the plane containing $\ff$. When $\pi\cap\bar{G}=\pi\cap\partial G=\ff$,
we call $\ff$ as a convex ``complete" face of $G$. If $\pi\cap\bar{G}\not=\pi\cap\partial G$, then $\ff$ is called a ``concave part" of $\partial G$.
In particular, when $\ff$ is a convex ``complete" face of some $G_r$, then $\ff$ is called a
concave ``local" face of $G$ (see Figure 4). A point $\vv\in\partial G$ is called a concave vertex of $G$ if $\vv$ is a common vertex of several concave ``local" faces of $G$;
similarly, an edge $\E\subset\partial G$ is called a concave edge of $G$ if $\E$ is a common edge of several concave ``local" faces of $\partial G$. For a
concave vertex $\vv$ (resp. concave edge $\E$) of $G$, the connected union of all concave ``local" faces that contain $\vv$ (resp. $\E$) as their common vertex (resp. common edge)
is called a concave ``complete" face of $G$.

\begin{center}
\includegraphics[width=6cm,height=6cm]{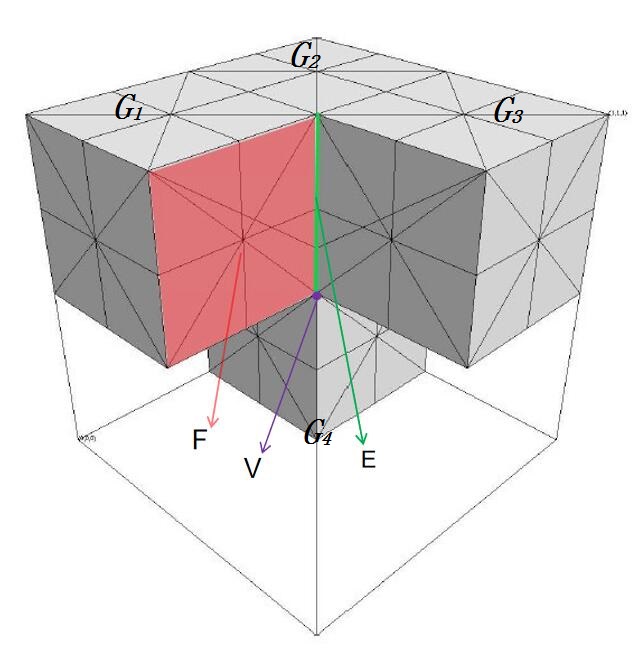}
\end{center}
\vskip 0.1in
Figure 4: $G$ denotes the union of four cubes $G_1,G_2,G_3,G_4$. The red square $F$ is a convex ``complete" face of $G_1$ but it is a concave ``local" face of $G$, the purple point $V$
is the unique concave vertex of $G$ and the green edge $E$ is a concave edge of $G$.\\
\vskip 0.1in

Let all the convex polyhedra $\{G_r\}_{r=1}^{m_0}$ be divided into two sets $\Sigma_1$ and $\Sigma_2$. For convenience, set
$$ \Lambda_k=\{r:~G_r\in\Sigma_k\}\quad\mbox{and}\quad\Gamma_k=\bigcup_{r\in\Lambda_k}\partial{G}_r\quad\quad(k=1,2). $$

Let $\Gamma$ be a connected ``Lipschitz" union of some complete faces.
We assume that the sets $\Sigma_1$ and $\Sigma_2$ satisfy the conditions: (i) the set $\Sigma_1$ consists of polyhedrons
that do not intersect each other, and the intersection of $\Gamma$ with each polyhedron in $\Sigma_1$ is either an empty set or a connected union of faces of the polyhedron;
(ii) the intersection of each polyhedron in $\Sigma_2$ with $\Gamma\cup\Gamma_{1}$ is a connected union of faces of the polyhedron;
(iii) the intersection of any polyhedron $G_l$ in $\Sigma_2$ with another polyhedron in $\Sigma_2$
is a subset of $G_l\cap(\Gamma\cup\Gamma_1)$ (which is a connected union of faces of $G_l$) but cannot be one face.

The following result gives a slightly weak $L^2$
stability of the regular decomposition on a class of non-convex polyhedra.

\begin{theorem}\label{regular decomposition-new} Let $G$ be a usual non-convex polyhedron, and $\Gamma$ be a connected ``Lipschitz" union of
some complete faces of $G$, with $G$ and $\Gamma$ satisfying the above assumptions.
Assume that $\v_h\in V_h(G)$
satisfies $\v_h\times\n=\0$ on $\Gamma$. Then there exist
$p_h\in Z_h(G)\cap H^1_{\Gamma}(G)$, $\w_h\in (Z_h(G)\cap H^1_{\Gamma}(G))^3$ and ${\bf
R}_h\in V_h(G)\cap H_{\Gamma}(\c;~G)$ such that \ee \v_h=\nabla
p_h+\r_h\w_h+{\bf R}_h,\label{decom4.new1}\e with the following
estimates \ee \l{estim4.new1}
\|\w_h\|_{1,G}+h^{-1}\|{\bf R}_h\|_{0,G}\stl\log(1/h)\|\c\,\v_h\|_{0,G}\e
and \ee
\|\w_h\|_{0,G}+\|p_h\|_{1,G}\stl\log(1/h)\|\v_h\|_{0,G}. \label{estim4.new2}\e
In particular, when $\Gamma$ contains all the ``concave parts" of $\partial G$, the above logarithmic factor
can be dropped.
\end{theorem}
{\it Proof}. When $\Gamma$ contains all the ``concave parts" of $\partial G$, the results directly follow from Lemma 4.2 and Example 4.1.
We only need to consider the case that $\Gamma$ does not contain the concave parts of $\partial G$. In order to make the ideas be understood more easily,
we first assume that $G$ is a union of three cubes: $G=D_1\cup D_2\cup D_3$ with $D_1=[0,\12]^3$, $D_2=[\12,1]\times[0,\12]^2$ and
$D_3=[0,\12]\times[\12,1]\times[0,\12]$ (see Figure 5). For this case, we can define $\Sigma_1=\{D_1\}$ and $\Sigma_2=\{D_2,~D_3\}$.
\begin{center}

\includegraphics[width=6cm,height=6cm]{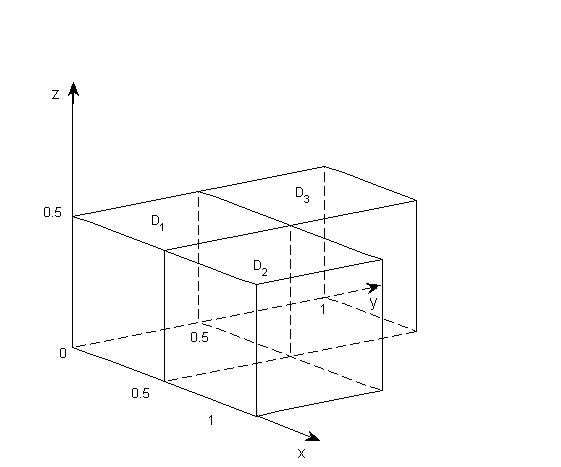}

\centerline{Figure 5: A non-convex polyhedron composed of three cubes}
\end{center}

We divide the proof into two steps.

{\bf Step 1}. Build the desired decomposition.

We first build a decomposition of $\v_h$ on $D_1$ by Lemma \ref{regular decomposition-3.new}.
The function $\v_{h,1}=\v_h|_{D_1}$ admits
the decomposition \ee \v_{h,1}=\nabla p_{h,1}+\r_h\w_{h,1}+{\bf
R}_{h,1}~~\mbox{on}~D_1,\l{decom4.new2}
\e
with $p_{h,1}\in
Z_h(D_1)$, $\w_{h,1}\in (Z_h(D_1))^3$ and ${\bf R}_{h,1}\in
V_h(D_1)$, which satisfy $p_{h,1}=0$, $\w_{h,1}=\0$ and ${\bf
R}_{h,1}\times\n=0$ on $\partial D_1\cap\Gamma$ (when
$\partial D_1\cap\Gamma=\emptyset$, we can require that $p_{h,1}$ has the zero average value on $D_1$). Moreover, since $D_1$ is convex (which implies the $L^2$ stability of the decomposition), we have
\ee \l{estim4.new3}
\|\w_{h,1}\|_{1,D_1}\stl\|\c\,\v_{h,1}\|_{0,D_1},\quad
\|\w_{h,1}\|_{0,D_1}+\|p_{h,1}\|_{1,D_1}\stl\|\v_{h,1}\|_{0,D_1}\e and \ee h^{-1}\|{\bf
R}_{h,1}\|_{0,D_1}\stl\|\c\,\v_{h,1}\|_{0,D_1}.
\label{estim4.new4}\e

Secondly, we extend $\w_{h,1}$ and $p_{h,1}$ into $D_2$ and $D_3$ in a special manner such that
some stability can be satisfied.

For $k=2,3$, set $\ff_{1k}=\partial D_1\cap\partial D_k$. By Lemma 4.5, there exists a function $\w^{\ff_{1k}}_{h,1}\in (Z_h(D_1))^3$
such that $\w^{\ff_{1k}}_{h,1}$ equals $\w_{1k}$ at all the nodes in the interior of $\ff_{1k}$ and $\w^{\ff_{1k}}_{h,1}=\0$ on $\partial D_1\backslash\ff_{1k}$ ($k=2,3$).
By the extension theorem and the Scott-Zhang interpolation \cite{ScottZ1990}
, we can show  (refer to the proof of Lemma 4.5 in \cite{Klaw2008}) there exists
an extension $\tilde{\w}^{\ff_{1k}}_{h,1}\in (Z_h(G))^3$ such that $\tilde{\w}^{\ff_{1k}}_{h,1}=\w^{\ff_{1k}}_{h,1}$ on $D_1$, $\tilde{\w}^{\ff_{1k}}_{h,1}$
vanishes on $\partial D_k\backslash\ff_{1k}$ ($k=2,3$) and satisfies
\ee
\|\tilde{\w}^{\ff_{1k}}_{h,1}\|_{1,D_k}\stl\|\w^{\ff_{1k}}_{h,1}\|_{1,D_1}~~~\mbox{and}~~~\|\tilde{\w}^{\ff_{1k}}_{h,1}\|_{0,D_k}\stl\|\w^{\ff_{1k}}_{h,1}\|_{0,D_1}~~~~(k=2,3).
\label{stab:3.0001new}
\e
Here we have used $L^2$ approximation of the quasi-interpolation and the inverse inequality to derive the $L^2$ stability (refer to (\ref{stab-lemma3.2-new})).
Set $\ff^{\partial}=\partial\ff_{12}\cup\partial\ff_{13}$, and let
$\tilde{\w}^{\partial}_{h,1}\in (Z_h(G))^3$ denote the natural zero
extension of $\w_{h,1}|_{\ff^\partial}$. Define $\tilde{\w}_{h,1}$ as follows: \
$$ \tilde{\w}_{h,1}=\w_{h,1}~~\mbox{on}~~D_1;~~
\tilde{\w}_{h,1}=\tilde{\w}^{\ff_{1k}}_{h,1}+\tilde{\w}^{\partial}_{h,1}|_{D_k}~~~\mbox{on}~~D_k~~(k=2,3). $$
It is easy to see that $\tilde{\w}_{h,1}\in (Z_h(G))^3$.

We define an extension $\tilde{p}_{h,1}\in Z_h(G)$ as follows: $\tilde{p}_{h,1}=p_{h,1}$ on $\bar{D}_1$;
$\tilde{p}_{h,1}$ vanishes at all the nodes in $\partial D_k\backslash\bar{\ff}_{1k}$ ($k=2,3$); $\tilde{p}_{h,1}$ is discrete harmonic
in $D_k$ ($k=2,3$). Let $\tilde{\bf R}_{h,1}\in V_h(G)$
be the natural zero extension of ${\bf R}_{h,1}$. For $k=2,3$, we define
\ee
\v^{\ast}_{h,k}=\v_h|_{D_k}-(\nabla \tilde{p}_{h,1}+\r_h\tilde{\w}_{h,1}+\tilde{\bf R}_{h,1})|_{D_k}~~\mbox{on}~~D_k. \label{decom4.new3}
\e
It is easy to see that $\v^{\ast}_{h,k}\times\n=\0$ on $\bar{\ff}_{1k}\cup(\partial D_k\cap\Gamma)$ ($k=2,3$).

Now we build the desired decomposition based on a regular decomposition of the function $\v^{\ast}_{h,k}$ ($k=2,3$) defined above.

It follows by Lemma \ref{regular decomposition-3.new} that the function $\v^{\ast}_{h,k}$ admits a decomposition
\ee \v^{\ast}_{h,k}=\nabla
p^{\ast}_{h,k}+\r_h\w^{\ast}_{h,k}+{\bf
R}^{\ast}_{h,k}~~\mbox{on}~D_k~~(k=2,3)\l{decom4.new5}\e
with $p^{\ast}_{h,k}\in Z_h(D_k)$,
$\w^{\ast}_{h,k}\in (Z_h(D_k))^3$ and ${\bf
R}^{\ast}_{h,k}\in V_h(D_k)$ ($k=2,3$), which satisfy
$p^{\ast}_{h,k}=0$, $\w^{\ast}_{h,k}=\0$ and
${\bf R}^{\ast}_{h,k}\times\n=0$ on $\bar{\ff}_{1k}\cup(\partial
D_k\cap\Gamma)$ ($k=2,3$). Moreover, notice that $D_k$ is a convex polyhedron (which implies the $L^2$ stability of the decomposition), for $k=2,3$ we have
\ee \l{estim4.new5} \|\w^{\ast}_{h,k}\|_{1,
D_k}\stl\|\c\,\v^{\ast}_{h,k}\|_{0,D_k},\quad
\|\w^{\ast}_{h,k}\|_{0,D_k}+\|p^{\ast}_{h,k}\|_{1,D_k}\stl\|\v^{\ast}_{h,k}\|_{0,D_k}\e
and \ee h^{-1}\|{\bf
R}^{\ast}_{h,k}\|_{0,D_k}\stl\|\c\,\v^{\ast}_{h,k}\|_{0,D_k}.
\label{estim4.new6}\e

Since $p^{\ast}_{h,k}$, $\w^{\ast}_{h,k}$ and
${\bf R}^{\ast}_{h,k}$ have the zero degrees of freedom on $\bar{\ff}_{1k}$, we can naturally extend them onto $G$ by zero.
We denote the resulting zero extentions by $\tilde{p}^{\ast}_{h,k}$, $\tilde{\w}^{\ast}_{h,k}$ and
$\tilde{\bf R}^{\ast}_{h,k}$. Define
\[p_h=\tilde{p}_{h,1}+\sum\limits_{k=2}^3\tilde{p}^{\ast}_{h,k},~~\w_h=\tilde{\w}_{h,1}+\sum\limits_{k=2}^3\tilde{\w}^{\ast}_{h,k}
~~\mbox{and}~~{\bf R}_h=\tilde{\bf R}_{h,1}+\sum\limits_{k=2}^3\tilde{\bf
R}^{\ast}_{h,k}.\]
It is easy to see that $p_h$, $\w_h$ and ${\bf R}_h$ have the zero degrees of freedom on $\Gamma$.
Using the local decompositions (\ref{decom4.new2}) and (\ref{decom4.new5}), together with the relation (\ref{decom4.new3}), we get
the global decomposition of $\v_h$
\ee \v_h=\nabla p_h+\r_h\w_h+{\bf R}_h.
\l{decom4.new6}\e

{\bf Step 2}. Derive the stability estimates.

From the definition of $\w_h$, we have
\ee
\|\w_h\|_{1,G}\stl\|\tilde{\w}_{h,1}\|_{1,G}
+\sum\limits_{k=2}^3\|\tilde{\w}^{\ast}_{h,k}\|_{1,D_k}.\label{ineq4.new1}
\e
For $k=2,3$, by (\ref{estim4.new5}) and (\ref{decom4.new3}) we can deduce that
$$\|\tilde{\w}^{\ast}_{h,k}\|_{1,D_k}\stl\|\c\,\v^{\ast}_{h,k}\|_{0,D_k}\stl\|\c\,\v_h\|_{0,D_k}+
\|\c(\r_h\tilde{\w}_{h,1})\|_{0,D_k}+\|\c\,\tilde{\bf R}_{h,1}\|_{0,D_k}.$$
Applying Proposition 4.1 and the inverse inequality, we further get
\beq
\|\tilde{\w}^{\ast}_{h,k}\|_{1,D_k}&\stl&\|\c\,\v_h\|_{0,D_k}+\|\c\,\tilde{\w}_{h,1}\|_{0,D_k}+h^{-1}\|\tilde{\bf R}_{h,1}\|_{0,D_k}\cr
&\stl&\|\c\,\v_h\|_{0,D_k}+\|\tilde{\w}_{h,1}\|_{1,D_k}+h^{-1}\|{\bf R}_{h,1}\|_{0,D_1}.\label{ineq4.new2}
\eq
Here we have used the relation $\|\tilde{\bf R}_{h,1}\|_{0,D_k}\stl\|{\bf R}_{h,1}\|_{0,D_1}$, which can be verified directly by the definition of $\tilde{\bf R}_{h,1}$.
Substituting (\ref{ineq4.new2}) into (\ref{ineq4.new1}), and using (\ref{estim4.new3})-(\ref{estim4.new4}) yields
\ee
\|\w_h\|_{1,G}\stl\|\c\,\v_h\|_{0,G}+\sum\limits_{k=2}^3\|\tilde{\w}_{h,1}\|_{1,D_k}.\label{ineq4.new3}
\e
Similarly, we can show
\ee
\|\w_h\|_{0,G}\stl\|\v_h\|_{0,G}+\sum\limits_{k=2}^3\|\tilde{\w}_{h,1}\|_{0,D_k}.\label{ineq:4.0003new}
\e

It suffices to estimate $\|\tilde{\w}_{h,1}\|_{1,D_k}$ and $\|\tilde{\w}_{h,1}\|_{0,D_k}$ ($k=2,3$). It follows by Lemma 4.5 that
$$
\|\w_{h,1}^{\ff_{1k}}\|_{1,D_1}\stl \log(1/h)\|\w_{h,1}\|_{1,
D_1}\quad (k=2,3).$$
Combining this inequality with (\ref{stab:3.0001new}) leads to
\ee
 \|\tilde{\w}_{h,1}^{\ff_{1k}}\|_{1,D_k}\stl \log(1/h)\|\w_{h,1}\|_{1,
D_1}.
\label{ineq:4.0001new}
\e
On the other hand, from Lemma 4.3, we have
$$\|\tilde{\w}^{\partial}_{h,1}\|_{1,D_k}\stl\|\w_{h,1}\|_{0,\ff^{\partial}}\stl\log^{\12}(1/h)\|\w_{h,1}\|_{1, D_1}~~~(k=2,3).$$
By the definition of $\tilde{\w}_{h,1}$, together with (\ref{ineq:4.0001new}) and the above estimate, we deduce that
$$ \|\tilde{\w}_{h,1}\|_{1,D_k}\stl\log(1/h)\|\w_{h,1}\|_{1, D_1}~~~(k=2,3).$$
Plugging this in (\ref{ineq4.new3}) gives the first estimate of the inequality (\ref{estim4.new1}).

By Lemma 4.5 and the definition of $\tilde{\w}^{\partial}_{h,1}$, we get \ee \|\w_{h,1}^{\ff_{1k}}\|_{0,D_1}\stl
\|\w_{h,1}\|_{0, D_1}~~~\mbox{and}~~~\|\tilde{\w}^{\partial}_{h,1}\|_{0,D_k}\stl\|\w_{h,1}\|_{0, D_1}~~~(k=2,3).\label{ineq:4.0002new} \e
This, together with the second inequality in (\ref{stab:3.0001new}), leads to
\ee
\|\tilde{\w}_{h,1}\|_{0,D_k}\stl\|\w_{h,1}\|_{0, D_1}~~~(k=2,3). \label{ineq:4.0005new}
\e
Substituting this into (\ref{ineq:4.0003new}) yields
\ee
\|\w_h\|_{0,G}\stl\|\v_h\|_{0,G}.\label{ineq:4.0002new}
\e

In the following we estimate $\|p_h\|_{1,G}$. If suffices to consider $\|\tilde{p}_{h,1}\|_{1,D_k}$ for $k=2,3$.
Since $\tilde{p}_{h,1}$ is discrete harmonic in $D_k$, we have
\ee
\|\tilde{p}_{h,1}\|_{1,D_k}\stl\|\tilde{p}_{h,1}\|_{{1\over 2},\partial D_k}~~~(k=2,3).
\label{ineq4.new4}
\e
Define the interpolation operators $I^0_{\ff_{1k}}$ and $I^0_{\partial\ff_{1k}}$ as follows:
for $\psi_h\in Z_h(G)$, the function $I^0_{\ff_{1k}}\psi_h$ (resp. $I^0_{\partial\ff_{1k}}\psi_h$) equals $\psi_h$ at the nodes in the interior of $\ff_{1k}$ (resp. on $\partial\ff_{1k}$)
and vanishes at all the other nodes on $\partial D_k$ (resp. at all the nodes not on $\partial\ff_{1k}$).
From the definition of $\tilde{p}_{h,1}$, we have $\tilde{p}_{h,1}=I^0_{\ff_{1k}}\tilde{p}_{h,1}+I^0_{\partial\ff_{1k}}\tilde{p}_{h,1}$ on $\partial D_k$. Then, it follows by (\ref{ineq4.new4}) that
\beqx
\|\tilde{p}_{h,1}\|_{1,D_k}&\stl&\|I^0_{\ff_{1k}}\tilde{p}_{h,1}\|_{{1\over 2}, \partial D_k}+\|I^0_{\partial\ff_{1k}}\tilde{p}_{h,1}\|_{{1\over 2}, \partial D_k}\cr
&\stl&\|I^0_{\ff_{1k}}\tilde{p}_{h,1}\|_{H^{{1\over 2}}_{00}(\ff_{1k})}+\|\tilde{p}_{h,1}\|_{0,\partial\ff_{1k}}\cr
&=&\|I^0_{\ff_{1k}} p_{h,1}\|_{H^{{1\over 2}}_{00}(\ff_{1k})}+\|p_{h,1}\|_{0,\partial\ff_{1k}}.
\eqx
Therefore, by using Lemma 4.4 and Lemma 4.3, we further obtain
$$ \|\tilde{p}_{h,1}\|_{1,D_k}\stl\log(1/h)\|p_{h,1}\|_{1, D_1}~~~(k=2,3).$$
Then, as in the estimate for $\|\w_h\|_{1,G}$ (but (\ref{ineq:4.0005new}) needs to be used), we get
$$ \|p_h\|_{1,G}\stl\log(1/h)\|\v_h\|_{0,G}. $$
Combining this with (\ref{ineq:4.0002new}) gives the inequality (\ref{estim4.new2}). Moreover, we can similarly derive
the second estimate in the inequality (\ref{estim4.new1}) by (\ref{estim4.new4}) and (\ref{estim4.new6}).

Next we extend the above ideas to the general case. Due to the condition (i), we can first independently define a regular decomposition of $\v_h$ on every convex polyhedron in $\Sigma_1$ and extend the functions
defined by the regular decompositions to all the polyhedra in $\Sigma_2$ in the same manner considered above. Then, by the conditions (ii) and (iii), we separately build regular decompositions
of the residual functions on all the polyhedra in $\Sigma_2$. Finally we sum all the functions generated by these regular decompositions to get the desired regular decomposition.
The condition (ii) makes Lemma \ref{regular decomposition-3.new} can be used for the residual functions on the polyhedra in $\Sigma_2$. Moreover, the condition (iii) makes two possibly neighboring polyhedra (which have a common edge or a common vertex) in $\Sigma_2$ can be separately handled since the functions defined by regular decompositions for the residual functions vanish on the intersection of the two polyhedra.
\hfill
$\Box$

\begin{remark} The construction of the decomposition
(\ref{decom4.new6}) is a bit technical. The main difficulty comes
from the definition of the vector-valued function $\w_h$, which must
satisfy $L^2$ stability. A natural idea is to extend $\w_{h,1}$
onto $G$ such that the extension is discrete harmonic in $D_2$ and $D_3$,
but the resulting extension may not satisfy the $L^2$ stability (\ref{ineq:4.0005new}).
\end{remark}

\section{Proof of Theorem \ref{teor:helm1}}
\setcounter{equation}{0}
In this section we are devoted to the proof of Theorem \ref{teor:helm1}. Throughout this section we always assume that $G$ is a usual (simply-connected) polyhedron,
and a face means a complete closed face of $G$. For convenience, we use the same notation $\r_h$ to denote the standard edge
element interpolation operator (if a connected union of some faces and edges of $G$ is involved) or the edge
element projector ${\bf R}^1_D$ introduced in \cite{Hiptmair-Pechstein2017} (if a non-connected union of some faces and edges of $G$ is involved). Of course, we can always
choose $\r_h$ as the projector ${\bf R}^1_D$.

 Let $\m{div}_{\tau}$ be the tangential divergence defined in
\cite{Alo1}, which was called surface $curl$ in \cite{Tos}. For
$\v_h\in V_h(G)$, we have
$curl_S\v_h=\m{div}_{\tau}(\n\times\v_h)=({\bf curl}~\v_h)\cdot\n$;
see \cite{Alo1}. For the purpose of clarity, we directly use the
notation $({\bf curl}~\v_h)\cdot\n$ in the rest of this paper.

%

Let ${\bf E}_h: V_h(\partial G)\rig V_h(G)$ be the discrete ${\bf curl}$-harmonic extension operator defined
in \cite{Ains2016}. This following result can be seen from the proof (see (4.2)) of Theorem 2.3 in \cite{Ains2016}. \\
{\bf Proposition 5.1} (stability of the ${\bf curl}$-harmonic
extension) ~For any $\v_h\in V_h(G)$, we have\ee \|\c~({\bf
E}_h(\v_h\times\n)|_{\partial G})\|_{0, G}\stl \|({\bf
curl}~\v_h)\cdot\n\|_{-{1\over 2},\partial G} \l{exten} \e
and
\ee\|{\bf
E}_h(\v_h\times\n)|_{\partial G}\|_{0,~G}\stl\|({\bf
curl}~\v_h)\cdot\n\|_{-{1\over 2},\partial G}+\|\v_h\times\n\|_{-{1\over 2}, \partial G}.\l{exten-new}
\e
\hfill $\Box$

 In the rest of this paper, we will repeatedly use $H^{-{1\over 2}}$-norm on a part of $\partial G$. For a connected ``Lipschitz" union $\Gamma$
 of some faces of $G$, let $H^{-{1\over 2}}(\Gamma)$ denote the dual space of $H^{{1\over 2}}(\Gamma)$ and define the corresponding norm as
 $$ \|\mu\|_{-\12,\Gamma}=\sup_{v\in H^{{1\over 2}}(\Gamma)}{\langle\mu, v\rangle\over \|v\|_{\12,\Gamma}},\quad \mu\in H^{-{1\over 2}}(\Gamma).$$
Notice that $H^{-{1\over 2}}(\Gamma)$ is different from the dual space of $H^{{1\over 2}}_{00}(\Gamma)$, which is often denoted by $H^{-{1\over 2}}_{00}(\Gamma)$
(see, for example, Subsection 7.1 of \cite{Tos}).

The following result can be directly built by the known $H^{-{1\over 2}}$-stability of the face extension
(cf.~Lemma 6 in \cite{hull}, Lemma 6.7 in \cite{hz} and Lemma 7.2 in \cite{Tos})\\
{\bf Proposition 5.2} Let $\v_h\in V_h(G)$ and $\Gamma$ be a connected ``Lipschitz" union of some faces of $G$. Then
\ee
\|({\bf curl}~\v_h)\cdot\n\|_{-\12,\Gamma}\stl
\log(1/h)\|({\bf curl}~\v_h)\cdot\n\|_{-\12,\partial G}. \l{face-exten} \e

Based on the above results we can build a discrete regular decomposition preserving zero
tangential trace on the boundary of a face of $G$.

\begin{lemma}\l{lem:helm2} Let $\ff$ be a (closed) face of $G$. Assume that $\v_h\in
V_h(G)$ satisfies $\v_h\cdot\t_{\partial\ff}=0$ on
$\partial\ff$. Then there exist $p_h\in Z_h(G)$, $\w_h\in
(Z_h(G))^3$ and ${\bf R}_h\in V_h(G)$, which satisfy
$p_h=0$, $\w_h=\0$ and ${\bf R}_h\cdot\t_{\partial\ff}=0$ on
$\partial\ff$, such that
\ee
\l{eq:desired1} \v_h=\nabla p_h+\r_h\w_h+{\bf R}_h
\e
with the following estimates
\ee
\l{eq:desired2}
\|\w_h\|_{1,G}+h^{-1}\|{\bf R}_h\|_{0,G}\stl\log(1/h)\|\c\,\v_h\|_{0,G}
\e
and
\ee
\|\w_h\|_{0,G}+\|p_h\|_{1,G}\stl\log(1/h)\|\v\|_{\c,G}.\label{stab:4.new11} \e
\end{lemma}
\no{\it Proof}. We separate the proof into two steps.

\smallskip
{\bf Step 1}: Establish the desired decomposition.

Define $\ff_c=(\partial G\backslash\ff)\cup\partial\ff$, which is a connected ``Lipschitz" union of faces. From the assumptions, we have a decomposition
\ee \v_h=\v^{\ff}_{h}+\v^{\ff_c}_{h}, \label{decom-lemma4.1-01}
\e
where $\v^{\ff}_{h}\in V_h(G)\cap H_{\ff}(\c;~G)$ is discrete $\c-$harmonic on $G$, and $\v^{\ff_c}_{h}\in V_h(G)\cap H_{\ff_c}(\c;~G)$ (it may be not discrete $\c-$harmonic on $G$).
For convenience, we use $\Gamma$ to denote $\ff$ or $\ff_c$. By Lemma
\ref{regular decomposition-3.new}, the function $\v^{\Gamma}_{h}$ ($\Gamma=\ff, \ff_c$) admits a discrete regular
decomposition
\ee
\v^{\Gamma}_{h}=\r_h\Phi_{h}^{\Gamma}+\nabla p_{h}^{\Gamma}+{\bf
R}_{h}^{\Gamma}~~~~(\Gamma=\ff, \ff_c), \label{decom-lemma4.1-02}
\e
where $\Phi_{h}^{\Gamma}\in (Z_h(G)\cap H^1_{\Gamma}(G))^3$, $p_{h}^{\Gamma}\in Z_h(G)\cap H^1_{\Gamma}(G)$
and ${\bf R}_{h}^{\Gamma}\in V_h(G)\cap H_{\Gamma}(\c;~G)$ ($\Gamma=\ff, \ff_c$). Moreover, they possess the stability estimates
\ee \|\Phi^{\Gamma}_h\|_{1,G}+h^{-1}\|{\bf R}_{h}^{\Gamma}\|_{0, G}\stl\|\c\,\v_{h}^{\Gamma}\|_{0,G}~~~~(\Gamma=\ff, \ff_c)
\label{stab-lemma4.1-01}\e
and
\ee
\|\Phi_{h}^{\Gamma}\|_{0,G}+\|p_{h}^{\Gamma}\|_{1,G}\stl\|\v_{h}^{\Gamma}\|_{\c,G}~~~~(\Gamma=\ff, \ff_c)
.\l{stab-lemma4.1-02} \e

Define
$$ \w_h=\Phi_{h}^{\ff}+\Phi_{h}^{\ff_c},~~p_h=p_{h}^{\ff}+p_{h}^{\ff_c}~~\mbox{and}~~{\bf R}_h={\bf R}_{h}^{\ff}+{\bf R}_{h}^{\ff_c}.$$
Then, from (\ref{decom-lemma4.1-01}) and (\ref{decom-lemma4.1-02}), these functions satisfy the regular decomposition (\ref{eq:desired1}). Moreover,
the functions $\w_{h}$, $p_h$ and ${\bf R}_h$ have zero degrees of freedom on $\partial\ff$ since $\partial\ff=\ff\cap\ff_c$.

\smallskip
{\bf Step 2}:  Verify the desired estimate (\ref{eq:desired2}) for
the decomposition (\ref{eq:desired1}).

From the definition of $\w_h$ and the triangle inequality, we have
$$
\|\w_h\|_{1, G}\stl\|\Phi_{h}^{\ff}\|_{1,G}+\|\Phi_h^{\ff_c}\|_{1, G}.
$$
This, along with (\ref{stab-lemma4.1-01}), leads to
\ee\|\w_h\|_{1,~G}\stl\|\c\,\v_{h}^{\ff}\|_{0,G}+\|\c\,\v_{h}^{\ff_c}\|_{0, G}. \label{stab-lemma4.1-03}
\e
By the definition of $\v_{h}^{\ff}$ and using the stability (\ref{exten}) of $\c$-harmonic extension,
we have
 \ee
\|\c\,\v_{h}^{\ff}\|_{0,G}\stl\|({\bf curl}~\v_{h}^{\ff})\cdot\n\|_{-\12,\partial G}\stl\|({\bf
curl}~\v_h)\cdot\n\|_{-\12,\ff_c}.
\l{stab4.7}\e
Moreover, using the estimate (\ref{face-exten}) and the trace inequality for $H(\c)$,
we deduce
$$ \|({\bf curl}~\v_h)\cdot\n\|_{-\12,\ff_c}\stl
\log(1/h)\|({\bf curl}~\v_h)\cdot\n\|_{-\12,\partial G}
\stl\log(1/h)\|\c\,\v_h\|_{0,G}. $$
Substituting this into (\ref{stab4.7}) yields
$$ \|\c\,\v_{h}^{\ff}\|_{0,G}\stl\log(1/h)\|\c\,\v_h\|_{0,G}. $$
This, together with (\ref{decom-lemma4.1-01}), further gives
$$ \|\c\,\v_{h}^{\ff_c}\|_{0,G}\stl\log(1/h)\|\c\,\v_h\|_{0,G}. $$
Combining (\ref{stab-lemma4.1-03}) with the above two inequalities, we obtain the first estimate of the inequality (\ref{eq:desired2}). The second estimate in (\ref{eq:desired2}) and
the inequality (\ref{stab:4.new11}) can be derived similarly (the stability (\ref{exten-new}) and the trace inequality
need to be used).
\hfill $\Box$

\begin{remark}
The conclusions in Lemma \ref{lem:helm2} are also valid when $\ff$ is replaced by a {\it connected ``Lipschitz" union of some faces} of $G$ (the union must has one connected boundary, so $\ff$ cannot be $\partial G$ itself).
But, if $\ff$ is a non-connected union of some faces, then the semi-norm on the right side of (\ref{eq:desired2}) must be replaced by the complete norm
.
\end{remark}

\begin{remark} Can the $\c$ norm on the right side of (\ref{stab:4.new11}) be replaced by the $L^2$ norm when the semi-norm on the right side of (\ref{eq:desired2}) is replaced
by the complete norm (comparing Lemma \ref{regular decomposition} and Lemma \ref{regular decomposition-3.new})?
Unfortunately, we failed to get the positive answer (notice that the $\c$ norm on the right side of the trace inequality cannot be replaced by the $L^2$ norm).
\end{remark}

From the above proof, we can obtain the following result\\
{\bf Corollary 5.1}. Let $\Gamma$ be a non-connected union of $\Gamma_1,\cdots,\Gamma_J$, where each $\Gamma_j$ is {\it connected ``Lipschitz" union of some faces} of $G$,
and let $\ff$ be a (closed) face of $G$ satisfying $\ff\cap\Gamma=\emptyset$.
Assume that $\v_h\in V_h(G)$ satisfies $\v_h\cdot\t_{\partial\ff}=0$ on $\partial\ff$ and $\v_h\times\n=\0$ on $\Gamma$. Then there exist $p_h\in Z_h(G)$,
$\w_h\in (Z_h(G))^3$ and ${\bf R}_h\in V_h(G)$, which satisfy
$p_h=0$, $\w_h=\0$ on $\Gamma\cup\partial\ff$ and $\lambda_e({\bf R}_h)=0$ for any $e\subset
\Gamma\cup\partial\ff$, such that
\ee \l{decomp:4.new1} \v_h=\nabla
p_h+\r_h\w_h+{\bf R}_h.\e Moreover, we have the following estimates
\ee
\l{stab:4.new21}
\|\w_h\|_{1,G}+h^{-1}\|{\bf R}_h\|_{0,G}\stl\log(1/h)\|\v_h\|_{\c,G}
\e
and
\ee
\|\w_h\|_{0,G}+\|p_h\|_{1,G}\stl\log(1/h)\|\v_h\|_{\c,G}.\label{stab:4.new22} \e \\
If $J=1$ and $\ff\cap\Gamma$ is an edge of $\ff$, then the norm on the right side of (\ref{stab:4.new21}) can be replaced by the $\c$ semi-norm.\\
\no{\it Proof}. As in the proof of Lemma \ref{lem:helm2}, we set $\ff_c=(\partial G\backslash\ff)\cup\partial\ff$ and use Lemma \ref{regular decomposition} for $\ff_c$ and $\Gamma\cup\ff$
(since Assumption 4.1 is met for them), respectively.
If $J=1$ and $\ff\cap\Gamma$ is an edge of $\ff$, we can use Lemma \ref{regular decomposition-3.new} since $\ff\cup\Gamma$ is connected. \hfill$\Box$

By Lemma \ref{lem:helm2} we can further build a regular decomposition preserving zero
tangential component on one edge.

\begin{lemma}\l{lem:helm3} Let $\E$ be a (closed) edge of $G$,
and $\v_h$ be a finite element function in $V_h(G)$ such that
$\v_h\cdot\t_{\E}=0$ on $\E$. Then $\v_h$ admits a
decomposition
$$ \v_h=\nabla p_h+\r_h\w_h+{\bf R}_h $$
for some $p_h\in Z_h(G)$, $\w_h\in (Z_h(G))^3$ and
${\bf R}_h\in V_h(G)$, which satisfy $p_h=0$, $\w_h=\0$,
${\bf R}_h\cdot\t_{\E}=0$ on $\E$. Moreover, the following estimates hold
\ee \l{eq:separate1}
\|\w_h\|_{1,G}+h^{-1}\|{\bf R}_h\|_{0,G}\stl\log(1/h)\|\c~\v_h\|_{0,G}
\e
and
\ee
\|\w_h\|_{0,G}+\|p_h\|_{1,G}\stl\log(1/h)\|\v\|_{\c,G}. \label{stab:4.new14}\e
The
conclusion is also valid for the case when $\E$ is replaced by a
connected union of several edges of $G$.
\end{lemma}

\no {\it Proof}. We separate the proof into three steps.

\smallskip
{\bf Step 1}:  Establish an edge-related decomposition.

Let $\ff$ be a face containing the edge $\E$. We first consider
a decomposition of the tangential component $\v_h\cdot\t_{\partial\ff}$
of $\v_h$ on $\partial\ff$. For convenience, we write
$\E_c=\partial\ff\backslash\E$.
Let $s$ be the arc-length variable along $\E_c$, taking values from
$0$ to $l_0$, where $l_0$ is the total length of $\E_c$. In terms of
$s$, the function $\v_h\cdot\t_{\E_c}$ is piecewise linear on the interval
$[0, l_0]$, denoted by $\hat v(s)$. Then we define
$$
C_{\E} = \frac{1}{l_0}\int_0^{l_0} \hat v(s)\,ds\quad\mbox{and}\quad
\phi_{\E}(t)
=\int_0^{t} (\hat v(s)-C_{\E})ds\,,
\quad \forall\,t\in [0, l_0]\,.
$$
Clearly we see $\phi_{\E}(t)$ vanishes at $t=0$ and $l_0$. We can
extend $\phi_{\E}$ naturally by zero onto $\E$, then
into $\partial G$ and $G$ such that its extension
$\tilde{\phi}_{\E}\in Z_h(G)$ and vanishes at all of the nodes that do not in the interior of $\E_c$. In the following, we define an
extension $\tilde C_{\E}$ of $C_{\E}$ such that $\tilde C_{\E}$
belongs to $(Z_h(G))^3$ and vanishes on $\E$. Moreover, we require that $\tilde
C_{\E}$ satisfies $(\r_h\tilde{C}_{\E})\cdot\t_{\partial\ff}=C_{\E}$ on
$\E_c=\partial\ff\backslash\E$ and \ee \| \tilde C_{\E}\|_{1,\hat
\o}\stl |C_{\E}|.\label{stab4.8new} \e
Let ${\mathcal N}_h^G$ denote the set of
the nodes on $G$, and let ${\mathcal N}_h^{\E_c}$ denote the set of the nodes
in $\E_c$. Then the values of the vector-valued function $\tilde
C_{\E}$ at all of the nodes in ${\mathcal N}_h^G\backslash{\mathcal N}_h^{\E_c}$ are defined to be
zero. Moreover, the values of the vector-valued function $\tilde
C_{\E}$ at the nodes in ${\mathcal N}_h^{\E_c}$ are defined such that $\tilde C_{\E}$
is linear on each (coarse) edge on $\E_c$, and $\|\tilde
C_{\E}\|^2_{0,\E_c}$ reaches the minimal value under the constraint
$(\r_h\tilde C_{\E})\cdot\t_{\partial\ff}=C_{\E}$ on $\E_c$. Notice that the number of degrees of freedom of the function $\tilde
C_{\E}$, which equals three times the number of vertices in $\E_c$, is greater than the number of coarse edges contained in $\E_c$.
Then the minimization problem (with a quadric objective functional and compatible linear constraints) has a solution. In
particular, if $C_{\E}=0$, the desired vector-valued function
$\tilde C_{\E}=\0$. Since each edge on $E^c$ is of size $O(1)$, we
have
$$ \|\tilde C_{\E}\|_{0,\E_c}\stl |C_{\E}|.$$
Moreover, by the definition of $\tilde C_{\E}$ and the discrete
norms, we get
$$ \| \tilde C_{\E}\|_{1,\hat
\o}\stl \|\tilde C_{\E}\|_{0,\E_c}.$$ Thus the inequality
(\ref{stab4.8new}) indeed holds. Using the definitions of $\tilde{\phi}_{\E}$ and $\tilde{C}_{\E}$, and noting that $\v_h\cdot\t_{\E}=0$, one can verify that
(cf.\,\cite{Tos})
\ee \v_h\cdot\t_{\partial\ff} =(\nabla
\tilde{\phi}_{\E})\cdot\t_{\partial\ff}
+(\r_h\tilde{C}_{\E})\cdot\t_{\partial\ff}.\l{4.002} \e
Moreover, by the discrete norms and Lemma 7.8 of \cite{Tos}, we can deduce that
\ee
\|\tilde{\phi}_{\E}\|_{1,G}\stl\|\phi_{\E}\|_{0,\partial\ff}\stl
\log(1/h)(\|\v_h\|_{0,G}+\|\c\,\v_h\|_{0,G}).\l{stab4.8}\e
For the derivation of the second inequality in (\ref{stab4.8}), we have used the following relations (see \cite{Tos} for the details)
$$ \|\phi_{\E}\|_{0,\partial\ff}\stl\|\v_h\cdot\t_{\partial\ff}-C_{\E}\|_{H^{-1}(\partial\ff)}\stl\|\v_h\times\n\|_{-{1\over 2},\ff}+\|(\c~\v_h)\cdot\n\|_{-{1\over 2}, \ff}. $$
The logarithm in (5.20) essentially comes from the $H^{-{1\over 2}}$-stability of the face extension of finite element functions.

\smallskip
{\bf Step 2}: Construct the desired decomposition in Lemma \ref{lem:helm3}.

For this purpose, we set
\begin{equation}\label{eq:step21}
\hat \v_{h,\E}=\v_h-(\nabla \tilde{\phi}_{\E}+\r_h\tilde{C}_{\E}).
\end{equation}
From (\ref{4.002}), we know that $\hat\v_{h,\E}\cdot\t_{\partial\ff}=0$ on
${\partial\ff}$. For the function $\hat \v_{h,\E}$ defined in
(\ref{eq:step21}), by Lemma~\ref{lem:helm2} one can find
functions $\hat p_h\in Z_h(G)$, $\hat \w_h\in
(Z_h(G))^3$ and $\hat{\bf R}_h\in V_h(G)$, which vanish on ${\partial\ff}$, such that
$$
\hat \v_{h,\E}=\nabla \hat p_h+\r_h\hat \w_h+{\bf \hat R}_h
$$
with the following estimates \ee
\|\hat\w_h\|_{1,G}\stl\log(1/h)\|\c\,\hat{\v}_{h,\E}\|_{0,G},~~~\|\hat p_h\|_{1,G}\stl\log(1/h)\|\hat{\v}_{h,\E}\|_{0,G}
\label{stab:4.new15} \e
and
\ee
h^{-1}\|{\bf \hat R}_h\|_{0,\hat
\o}\stl\log(1/h)\|\c\,\hat{\v}_{h,\E}\|_{0,\hat \o}.\l{4.006} \e

Now by defining
$$
p_h=\tilde{\phi}_{\E}+\hat p_h,~~\w_h=\tilde{C}_{\E}+\hat
\w_h~~\mbox{and}~~{\bf R}_h={\bf \hat R}_h\,,
$$
we get the final decomposition \ee\l{eq:decomp1}
 \v_h=\nabla p_h+\r_h\w_h+{\bf R}_h
\e such that $p_h=0$ and $\w_h=0$ on ${\E}$.

\smallskip
{\bf Step 3}: Derive the desired estimate in Lemma~\ref{lem:helm3} for the
decomposition (\ref{eq:decomp1}).

By the definition of $C_{\E}$ and noting that $\v_h\cdot\t_{\E}=0$ on ${\E}$, which implies that
$\v_h\cdot\t_{\partial\ff}=0$ on ${\E}$, we have
$$ C_{\E}={1\over l_0}\int_{\t_{\E_c}}\v_h\cdot\t_{\E_c}ds={1\over l_0}\int_{\partial\ff}\v_h\cdot\t_{\partial\ff}ds. $$
Thus, by the Stokes' theorem on $\ff$, we deduce
$$ C_{\E}={1\over
l_0}\int_{\ff}({\bf curl}~\v_h)\cdot\n ds. $$
Then
\ee |C_{\E}|\stl|\int_{\ff}({\bf curl}~\v_h)\cdot\n ds|\stl\|({\bf
curl}~\v_h)\cdot\n\|_{-{{1\over 2}}, \ff}\cdot\|1\|_{\12,~\ff}.
\l{5.001new}\e
Using the estimate (\ref{face-exten}) again, we have
$$
\|({\bf curl}~\v_h)\cdot\n\|_{-{{1\over 2}}, \ff}\stl
\log(1/h)\|({\bf curl}~\v_h)\cdot\n\|_{-\12,~\partial G},
$$
 This,  along with
(\ref{5.001new}), leads to
$$ |C_{\E}|\stl\log(1/h)\|\c~\v_h\|_{0, G}. $$
By (\ref{stab4.8new}), we further obtain
$$ \|
\tilde C_{\E}\|_{1,\hat \o}\stl\log(1/h)\|\c~\v_h\|_{0,
G}. $$ By the definition of $\w_h$ and the above inequality, together with (\ref{stab:4.new15}),
(\ref{eq:step21}) and (\ref{stab4.8}), we deduce the first estimate in the inequality (\ref{eq:separate1})

In a similar way, we can prove the second estimate in
(\ref{eq:separate1}) and the inequality (\ref{stab:4.new14}) by (\ref{4.006}).

When $\E$ is replaced by a connected union of some edges, the proof is almost the same as that given above, where the face $\ff$ may be replaced by a connected ``Lipschitz" union
of several faces (the union has one connected boundary, so it cannot be $\partial G$ itself). \hfill$\Box$

\begin{remark} There is a key difference in the proof of the above lemma from that of Lemma 4.3 in \cite{HuShuZou2013}:
since the extension $\tilde C_{\E}$ in the above proof must belong
to the space $(Z_h(G))^3$, $\tilde C_{\E}$ cannot be defined
to be the natural zero extension of $C_{\E}$ as in
\cite{HuShuZou2013}.
\end{remark}

From the explanations in Example 4.1, we know that Lemma \ref{regular decomposition}  and Lemma \ref{regular decomposition-3.new} do not hold
when some $\Gamma_r$ is a connected ``non-Lipschitz" union of some faces.
The following result can be regarded as extensions of Lemma \ref{regular decomposition} and Lemma \ref{regular decomposition-3.new}.

\begin{lemma}\label{lem:face-new} Let $\Gamma$ be a connected union of some faces of $G$. Assume that $\v_h\in V_h(G)$ has zero tangential trace on $\Gamma$.
Then the function $\v_h$ has a decomposition
$$ \v_h=\nabla p_h+\r_h\w_h+{\bf R}_h$$
with some $p_h\in Z_h(G)$, $\w_h\in (Z_h(G))^3$ and $
{\bf R}_h\in V_h(G)$ satisfying $p_h=0$ and $\w_h=\0$ on $\Gamma$. Moreover, we have
\ee
\|\w_h\|_{1,G}+h^{-1}\|{\bf
R}_h\|_{0,G}\stl\log(1/h)\|\c\,\v_h\|_{0,G}\label{stab:4.new016}\e
and
\ee
\|\w_h\|_{0,G}+\|p_h\|_{1,G}\stl\log(1/h)\|\v_h\|_{\c,G}. \label{stab:4.new017}\e
If $\Gamma$ is a non-connected union of $\Gamma_1,\Gamma_2,\cdots,\Gamma_J$, with each $\Gamma_r$ being a connected union of some faces of $G$, then the semi-norm on the right side of
(\ref{stab:4.new016}) should be replaced by the complete $\c$ norm.
\end{lemma}
\no {\it Proof}. Here we need not assume that $\Gamma$ is a connected ``Lipschitz" union of some faces, i.e., $\Gamma$ may contain one isolated vertex: there exists one face on $\Gamma$ such that
the intersection of this face with any neighboring face on $\Gamma$ is a vertex in $\Gamma$.
For this case, Assumption 4.1 is not met and so Lemma \ref{regular decomposition-3.new} cannot be used directly. For simplicity of exposition,
we only focus on the simplest case that $\Gamma\subset\partial G$ is a union of two faces that just have a common vertex.
There is no essential difference in the proof for more general cases
(see the explanations given later).

Let $\Gamma=\ff_1\cup\ff_2$ with $\ff_1$ and $\ff_2$ being two closed faces of $G$, where the intersection of $\ff_1$ and $\ff_2$ is a vertex $\vv$ (i.e., $\vv=\ff_1\cap\ff_2$ is an isolated vertex).
We use $\E_i\subset\ff_i$ ($i=1,2$) to denote one of the closed edges containing $\vv$ as their endpoint, such that $\E=\E_1\cup\E_2$ is on the same side of $\Gamma$. Then
there exists a connected ``Lipschitz" union $\hat{\Gamma}$ of some faces such that
$\E\subset\partial\hat{\Gamma}$ ($\hat{\Gamma}$ may be just one face), see Figure 6.

Since $\E$ is connected and $\v_h$ has zero tangential components on $\E$ by the assumption, we can define $\tilde{C}_{\E}\in (Z_h(G))^3$ and $\tilde{\phi}_{\E}\in Z_h(G)$ as in Lemma \ref{lem:helm3} (choosing $\ff=\hat{\Gamma}$). Define
$$\hat \v_{h,\E}=\v_h-(\nabla \tilde{\phi}_{\E}+\r_h\tilde{C}_{\E}).$$
\begin{center}
\includegraphics[width=7.5cm,height=7cm]{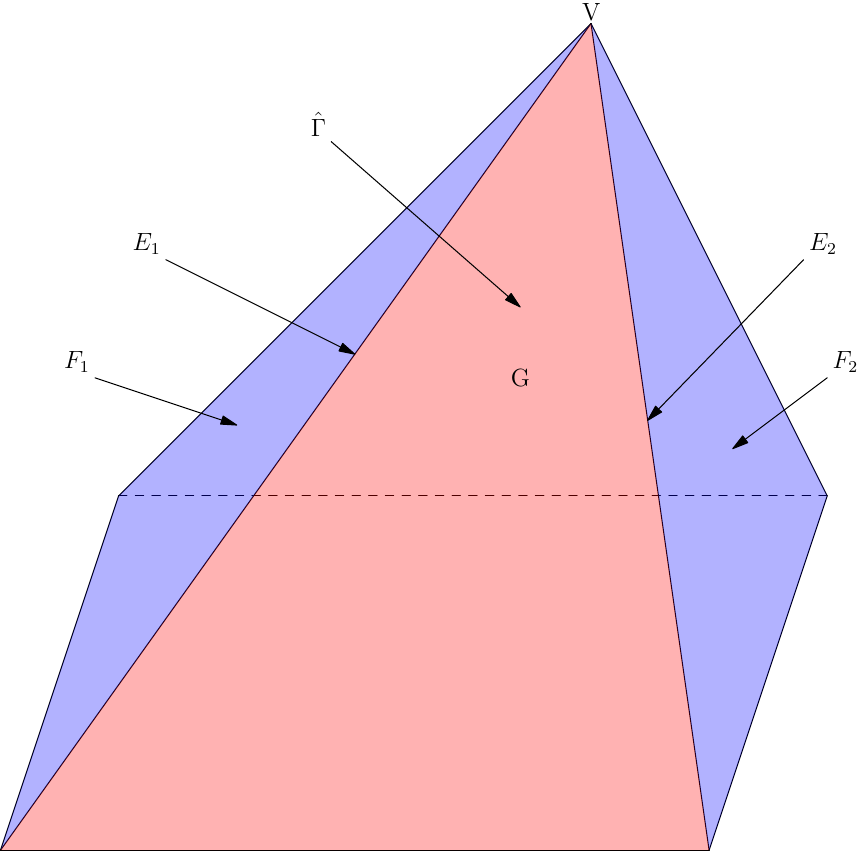}
\end{center}
Figure 6: The blue faces denote $\ff_1$ and $\ff_2$, which are two opposite lateral
faces of four pyramid domain $G$, the red lateral face denotes $\hat{\Gamma}$.
\vskip 0.2in

Then $\hat\v_{h,\E}\cdot\t_{\partial\hat{\Gamma}}=0$ on $\partial\hat{\Gamma}$. Set $\hat{\Gamma}_c=(\partial G\backslash\hat{\Gamma})\cup\partial\hat{\Gamma}$ and
$\tilde{\Gamma}=\Gamma\cup\hat{\Gamma}$. It is easy to see that both $\hat{\Gamma}_c$ and
$\tilde{\Gamma}$ are connected ``Lipschitz" unions of some faces (i.e., any vertex of them is not an isolated vertex). By the assumptions, we have a decomposition
$$ \hat{\v}_{h,\E}=\hat{\v}^{\hat{\Gamma}_c}_{h,\E}+\hat{\v}^{\tilde{\Gamma}}_{h,\E}, $$
where $\hat{\v}^{\hat{\Gamma}_c}_{h,\E}\in V_h(G)\cap H_{\hat{\Gamma}_c}(\c;~G)$ is discrete $\c-$harmonic on $G$, and $\hat{\v}^{\tilde{\Gamma}}_{h,\E}\in V_h(G)\cap H_{\tilde\Gamma}(\c;~G)$.
Then we can build discrete regular decompositions respectively for the two functions by Lemma \ref{regular decomposition-3.new}
since Assumption 4.1 is met for $\hat{\Gamma}_c$ and $\tilde{\Gamma}$. The resulting functions have zero traces on $\Gamma=\hat{\Gamma}_c\cap\tilde{\Gamma}$.
Finally, we can prove the desired results as in the proof of Lemma \ref{lem:helm3} and Lemma \ref{lem:helm2}.

If $\Gamma$ is a non-connected union of $\Gamma_1,\Gamma_2,\cdots,\Gamma_J$, with some $\Gamma_r$ being a union of two faces that just have a common vertex, we can
handle $\Gamma_r$ in the same manner as the one used above, and transform it into the cases that each considered part
is connected ``Lipschitz" union of some faces of $G$. Then we use Lemma \ref{regular decomposition} several times to get the desired results.
\hfill$\Box$

By Lemma \ref{lem:face-new} and the explanations in Example 4.1, the discrete results in Lemma \ref{regular decomposition-3.new} are essentially valid for any connected
union $\Gamma$ of faces of a polyhedron. With the help of Lemma \ref{lem:face-new}, we need not to emphasize the word ``Lipschitz" in {\it connected ``Lipschitz" union}
and simply say {\it connected union} in the following discussions.

\begin{lemma}\l{lem:helmFE} Let $\Gamma$ be a connected union of some faces and edges of $G$. Suppose that $\v_h$ satisfies {\bf Assumption 3.1}. Then $\v_h$ can be decomposed as
\ee \v_h=\nabla p_h+\r_h\w_h+{\bf R}_h\label{4.new0007}
\e
for some $p_h\in Z_h(G)$ and $\w_h\in (Z_h(G))^3$ and
${\bf R}_h\in V_h(G)$ such that $p_h$ and $\w_h$ vanish on $\Gamma$. Moreover, we have
 \ee \|\w_h\|_{1,G}+h^{-1}\|{\bf R}_h\|_{0,G}\stl\log(1/h)\|\c~\v_h\|_{0,G}
\label{4.new0008}\e
and
\ee
\|\w_h\|_{0,G}+\|p_h\|_{1,G}\stl\log(1/h)\|
\v_h\|_{\c,G}. \l{4.new0001}\e
The results are also valid when $\E$ is replaced by a connected union of several edges
of $G$. \end{lemma}
\no {\it Proof}.
Let $\E\subset\Gamma$ be an isolated edge (or a union of several edges) $\E$ of $\Gamma$, i.e., $\E$ is not an edge of any face in $\Gamma$. Since $\Gamma$ is connected, there exists a face $\ff\subset\Gamma$
(provided that $\Gamma$ indeed contains faces) such that the intersection $\E\cap\ff$
is an endpoint of $\E$. Let $\E'$ be an edge of the face $\ff$ such that $\E'\cap\E$ is just this endpoint, which implies that $\E'\cup\E$ is connected. Let $\hat{\Gamma}$ be a connected ``Lipschitz"
union of some faces of $G$ such that $\E\cup\E'\subset\partial\hat{\Gamma}$. By the assumptions, the function $\v_h$ has zero degrees of freedom on $\E\cup\E'$.

Without loss of generality, we assume that all of the faces in $\Gamma$ constitute a connected ``Lipschitz" union of some faces, otherwise, we need to add a union of some faces
as in the proof of Lemma \ref{lem:face-new}. Under this assumption, the set $\Gamma\cup\hat{\Gamma}$ is also a connected ``Lipschitz" union of some faces.
 Then the results in this lemma can be built as in the proof of Lemma \ref{lem:helm3} (replacing $\E$ by $\E\cup\E'$), by using Lemma \ref{regular decomposition-3.new}
on $\hat{\Gamma}_c=(\partial G\backslash\hat{\Gamma})\cup\partial\hat{\Gamma}$ and $\Gamma\cup\hat{\Gamma}$, respectively.
For this case, the norm on the right side of (\ref{4.new0008}) can be replaced by the $\c$ semi-norm.
\hfill $\Box$

\begin{lemma}\l{lem:helm5} Let $\E_1,\cdots,\E_n$ be (closed) edges of $G$, which satisfy $\E_i\cap\E_j=\emptyset$ for $i\not=j$. Assume that
$\v_h\in V_h(G)$ satisfies $\v_h\cdot\t_{\E_l}=0$ on each $\E_l$. Then $\v_h$ can be decomposed as
$$ \v_h=\nabla p_h+\r_h\w_h+{\bf R}_h$$
for some $p_h\in Z_h(G)$, $\w_h\in (Z_h(G))^3$ and
${\bf R}_h\in V_h(G)$ such that $p_h$ and $\w_h$ vanish on each edge $\E_l$. Moreover, we
have \ee \|\w_h\|_{1,G}\stl\log(1/h)\|\v_h\|_{\c,G},~~~\|p_h\|_{1,G}\stl\log(1/h)\|\v_h\|_{\c,G}\label{stab:4.new18}
\e
and
\ee
h^{-1}\|{\bf R}_h\|_{0,G}\stl\log(1/h)\|\v_h\|_{\c,G}. \l{4.000}\e
The results are also valid when $\E_l$ is replaced by a connected union of edges.
\end{lemma}
\no {\it Proof}. We first consider a simple case that, for each $\E_l$, there exists a face $\ff_l\subset\partial G$ such that $\ff_l$ contains $\E_l$ and the face $\ff_l$ does not adjoin the other edges.
For this case, regarding $\E_l$ as the edge $\E$ in the proof of Lemma \ref{lem:helm3} and almost repeating the proof process (but using Lemma \ref{regular decomposition} for $\Gamma=\cup_{l=1}^J
\ff_l$), we can obtain the desired results.

If the above condition is not met, the proof of this lemma is a bit technical. Without loss of generality, we assume that this condition is not satisfied for each $\E_l$ (An example
is that $\E_1,\E_2,\E_3,\E_4$ just are four parallel edges of a cube $G$, see Figure 7). This means that, for each edge $\E_l$,
any face containing $\E_l$ also contains another different edge $\E_{l'}$ (which does not intersect $\E_l$). In this situation, the above proof is not practical since the functions
$\tilde{\phi}_{\E}$ and $\tilde{C}_{\E}$ defined in the proof of Lemma \ref{lem:helm3} may not vanish on $\E_{l'}$.

\begin{center}
\includegraphics[width=7cm,height=6cm]{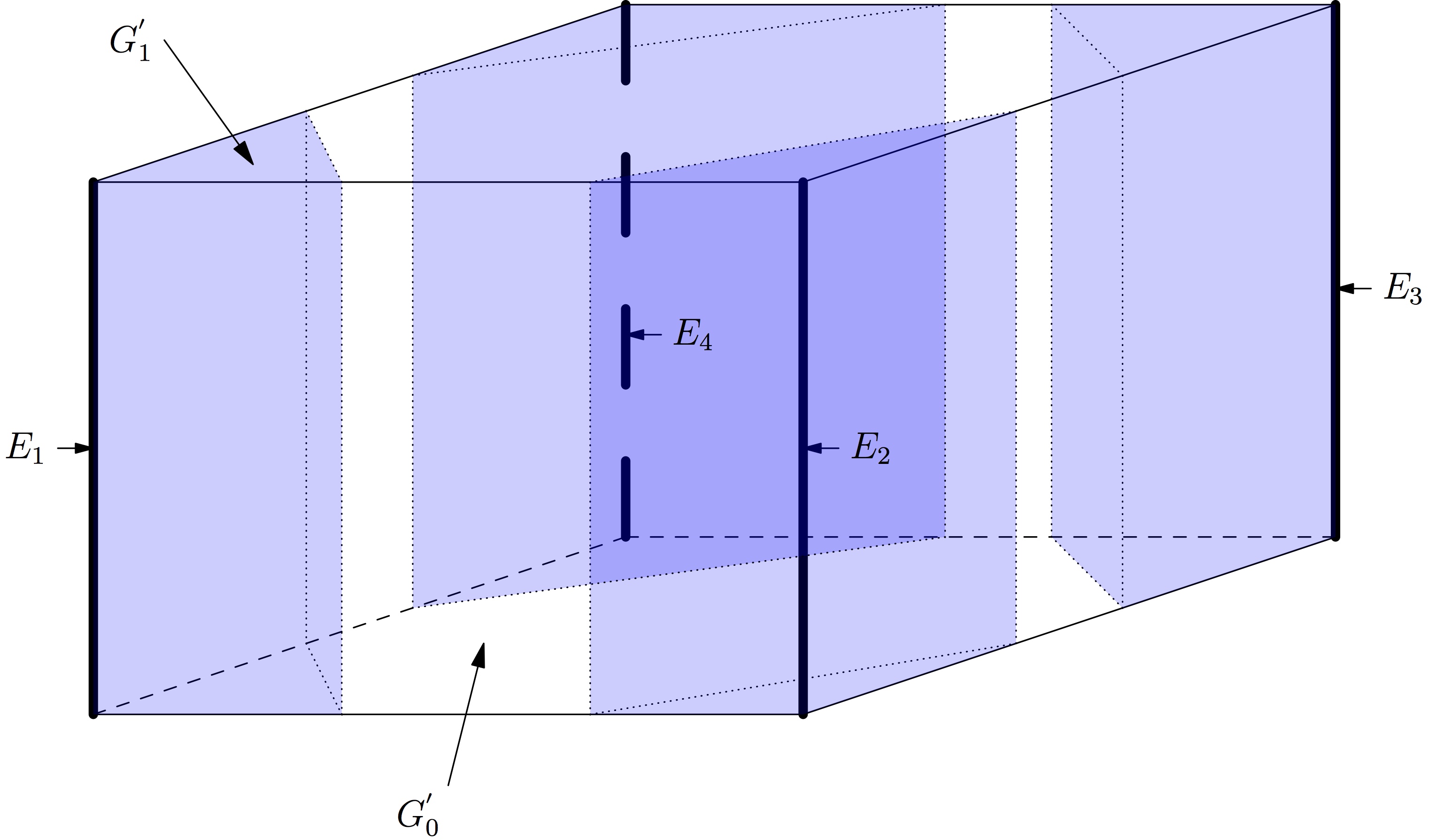}
\vskip 0.1in
\centerline{Figure 7: An example with four edges: the left shaded polyhedron denotes the}
\centerline{subdomain $G_1'$, the middle polyhedron denotes $G_0'$.}
\end{center}

Since the considered edges $\{\E_r\}$ do not intersect each other, we have $dist(\E_i,\E_j)=O(1)$ for $i\not=j$ (notice that $G$ is a usual polyhedron with a diameter $O(1)$).
Thus, we can decompose $G$ into a union of non-overlapping subdomains $G'_0, G'_1,\cdots,G'_m$ such that: (i) each subdomain $G'_l$
is a polyhedron with a diameter $O(1)$; (ii) $\bar{G}'_l\cap\bar{G}'_j=\emptyset$ for $j\not=l$ $(l,j\not=0$), and $G'_0$ just has a common face $\Gamma'_{0l}$ with each $G'_l$ ($l\not=0$);
(iii) for $l=1,\cdots,m$, the subdomain $G'_l$ contains $\E_l$ as one of its edges, but the subdomain $G'_0$ does not intersect any $\E_l$. In general we cannot require each subdomain
$G'_l$ to be a union of some elements. Because of this, for each $l$ we choose another subdomain $G_l$ generated from $G_l'$ (with a small perturbation only), where $G_l$ is a union
of all the elements $K$ satisfying $\mbox{meas}(K\cap G_l)\geq {1\over 2}\mbox{meas}(K)$, where $G_l$ can be also generated as in Example 4.2. It is clear that $G_l$ is not a usual polyhedron since
$\Gamma_{0l}=\bar{G}_0\cap \bar{G}_l$ is not a plane face yet. Fortunately, all the subdomains $\{G_l\}$ still constitute a union of $G$ and preserve the other properties of $\{G'_l\}$.

For each $\E_l$, we can use Lemma \ref{lem:helm3} to build a regular decomposition\footnote{Regarding $\E_l$ as $\E$ in Lemma \ref{lem:helm3}, we first get three finite element functions by the
regular decomposition of $\v_h$ on $G$, then we define $\w^{(l)}_h$, $p^{(l)}_h$ and ${\bf R}^{(l)}_h$ as the restrictions of the three functions on $G_l$, respectively.}

\ee
\v_h=\r_h\w^{(l)}_h+\nabla p^{(l)}_h+{\bf R}^{(l)}_h~~~\mbox{on}~~G_l \label{decom:4.newnew1}
\e
with $\w^{(l)}_h\in (Z_h(G_l))^3$ and $p^{(l)}_h\in Z_h(G_l)$, which vanish on $\E_l$ but may not vanish on the other edges. The decomposition is stable with a logarithmic factor.
Let $d_0$ be a given positive number independent of $h$. For $l=1,\cdots,m$,
we choose a ball $D_l$ containing $G_l$ such that $dist(\partial D_l, \partial G_l)\geq d_0$ and $D_l$ does not intersect any $G_j$ for $j\not=0,l$.
By the extension theorem, there exists an extension
$\tilde{\w}^{(l)}$ (resp. $\tilde{p}^{(l)}$) of $\w^{(l)}_h|_{G_l}$ (resp. $p^{(l)}_h|_{G_l}$) such that: (a) $\tilde{\w}^{(l)}\in (H^1(\mathbb{R}^3)^3$ (resp.
$\tilde{p}^{(l)}\in H^1(\mathbb{R}^3$); (b) $\tilde{\w}^{(l)}$ and $\tilde{p}^{(l)}$ vanish on the outside of $D_l$; (c) $\|\tilde{\w}^{(l)}\|_{1, G_0}\stl
\|\w^{(l)}_h\|_{1, G_l}$ and $\|\tilde{p}^{(l)}\|_{1, G_0}\stl\|p^{(l)}_h\|_{1, G_l}$. For each $l$, let $\tilde{\bf R}^{(l)}_h\in V_h(G)$ denote the standard zero extension
of ${\bf R}^{(l)}_h|_{G_l}$. Define
\ee
\tilde{\v}^{(0)}_h=\v_h-\sum\limits_{l=1}^m(\r_h\tilde{\w}^{(l)}+\nabla \tilde{p}^{(l)}+\tilde{\bf R}^{(l)}_h)\quad\mbox{on}~~G_0.\label{decom:4.newnew2}
\e

It is clear that $\tilde{\v}^{(0)}_h\in H(\c; G_0)$. Since $\tilde{\w}^{(l)},~\tilde{p}^{(l)}$ and $\tilde{\bf R}^{(l)}_h$ vanish on $\bar{G}_j$ for $j\not=0,l$, by (\ref{decom:4.newnew1})
we have $\tilde{\v}^{(0)}_h\times\n=\0$ on $\Gamma_{0l}$ for $l=1,\cdots,m$. From the explanations in Example 4.2, we know that Assumption 4.1 is satisfied for $G_0$. By Lemma 4.1, the function $\tilde{\v}^{(0)}_h$ admits a stable decomposition
\ee
\tilde{\v}^{(0)}_h=\tilde{\w}^{(0)}+\nabla \tilde{p}^{(0)}\quad\mbox{on}~~G_0 \label{decom:4.newnew3}
\e
for $\tilde{\w}^{(0)}\in (H^1(G_0))^3$ and $\tilde{p}^{(0)}\in H^1(G_0)$, with $\tilde{\w}^{(0)}$ and $\tilde{p}^{(0)}$ vanishing on $\Gamma_{0l}$ for $l=1,\cdots,m$. Combing (\ref{decom:4.newnew2})
and (\ref{decom:4.newnew3}), we have
$$\v_h=\tilde{\w}^{(0)}+\nabla \tilde{p}^{(0)}+\sum\limits_{l=1}^m(\r_h\tilde{\w}^{(l)}+\nabla \tilde{p}^{(l)}+\tilde{\bf R}^{(l)}_h)\quad\mbox{on}~~G_0. $$
Thus
\ee
\v_h=\r_h(\tilde{\w}^{(0)}+\sum\limits_{l=1}^m\tilde{\w}^{(l)})+\nabla p^{(0)}_h+\sum\limits_{l=1}^m\tilde{\bf R}^{(l)}_h\quad\mbox{on}~~G_0\label{decom:4.newnew4}
\e
with $p_h^{(0)}$ satisfying $\nabla p^{(0)}_h=\r_h\nabla(\tilde{p}^{(0)}+\sum\limits_{l=1}^m\tilde{p}^{(l)})$.

Let $\Pi_h: (H^1(G_0))^3\rightarrow (Z_h(G_0))^3$ denote the Scott-Zhang interpolation operator, which can preserve the values of a linear polynomial on some elements of the boundary $\partial G_0$.
Define
$$ \w^{(0)}_h=\Pi_h(\tilde{\w}^{(0)}+\sum\limits_{l=1}^m\tilde{\w}^{(l)})~~\mbox{and}~~{\bf R}^{(0)}_h=\r_h(I-\Pi_h)(\tilde{\w}^{(0)}+\sum\limits_{l=1}^m\tilde{\w}^{(l)})
+\sum\limits_{l=1}^m\tilde{\bf R}^{(l)}_h.$$
Then (\ref{decom:4.newnew4}) can be written as
\ee
\v_h=\r_h\w^{(0)}_h+\nabla p^{(0)}_h+{\bf R}^{(0)}_h\quad\mbox{on}~~G_0.\label{decom:4.newnew5}
\e
It is clear that $\w^{(0)}_h=\w^{(l)}_h$ and $p^{(0)}_h=p_h^{(l)}$ on $\Gamma_{0l}$ for $l=1,\cdots,m$, which implies that ${\bf R}^{(0)}_h\times\n={\bf R}^{(l)}_h$
on $\Gamma_{0l}$ for $l=1,\cdots,m$. Thus we naturally define $\w_h=\w^{(l)}_h$, $p_h=p_h^{(l)}$ and ${\bf R}_h={\bf R}^{(l)}_h$ on $G_l$ for $l=0,1,\cdots,m$, and we have
$\w_h\in (Z_h(G))^3$, $p_h\in Z_h(G)$ and ${\bf R}_h\in V_h(G)$, which have zero degrees of freedom on all of the edges $\E_l$.
It is easy to see from (\ref{decom:4.newnew1}) and
(\ref{decom:4.newnew5}) that the desired regular decomposition is valid for the defined functions. Besides, we can verify that the resulting regular decomposition is
also stable with a logarithmic factor.  \hfill $\Box$

\begin{remark} In the proofs of Lemma \ref{lem:helm2}-Lemma \ref{lem:helm5}, our main ideas are to transform the problem vanishing on an edge $\E$ into a
problem vanishing on a face $\ff$ of $G$ and then to use Lemma \ref{regular decomposition} or Lemma \ref{regular decomposition-3.new} to build a regular decomposition,
which can preserve zero trace on this face.
\end{remark}

Notice that, in all the previous Lemmas, a connected union of several edges  has no essential difference from an edge. Now we can prove Theorem \ref{teor:helm1} by
Lemma \ref{lem:helmFE} and Lemma \ref{lem:helm5}.

\vskip 0.2in

\no{\it Proof of Theorem \ref{teor:helm1}}.
We decompose $\Gamma$ into a  union of $\Gamma_1,\Gamma_2,\cdots,\Gamma_J$, where each $\Gamma_r$ is a connected union of some faces and edges of $G$ and
any two different $\Gamma_i$ and $\Gamma_j$ do not intersect (i.e., $\Gamma_i\cap\Gamma_j=\emptyset$).

Since $\Gamma_i\cap\Gamma_j=\emptyset$ and $G$ is a usual polyhedron with a diameter $O(1)$, we have $dist(\Gamma_i,\Gamma_j)=O(1)$ for $i\not=j$.
As in the second part in the proof of Lemma \ref{lem:helm5}, we can decompose $G$ into a union of $G_0,G_1,G_2,\cdots,G_J$ such that: (i) each subdomain $G_l$
has the size $O(1)$ and is a union of some elements; (ii) $\bar{G}_l\cap\bar{G}_j=\emptyset$ for $j\not=l$ $(l,j\not=0$), and $G_0$ just has a common face $\Gamma_{0l}$ with each $G_l$ ($l\not=0$);
(iii) for $l=1,\cdots,J$, $\Gamma_l\subset \partial G_l$, but the subdomain $G_0$ does not intersect any $\Gamma_l$. Here $G_l$ may not be a usual polyhedron yet.
For $l=1,\cdots,J$, we can build a stable discrete regular decomposition of $\v_h|_{G_l}$ as in Lemma \ref{lem:helmFE} (detailed explanations can refer to the footnote in the proof of Lemma \ref{lem:helm5})
$$ \v_h|_{G_l}=\r_h\w^{(l)}_h+\nabla p^{(l)}_h+{\bf R}^{(l)}_h $$
with $p^{(l)}_h\in Z_h(G_l)$ and $\w^{(l)}_h\in(Z_h(G_l))^3$ vanishing on $\Gamma_l$. We stably extend $\w^{(l)}_h$, $p^{(l)}_h$ and ${\bf R}^{(l)}_h$ into $\mathbb{R}^3$ (see the proof of Lemma \ref{lem:helm5}), and
define
$$
\tilde{\v}^{(0)}_h=\v_h-\sum\limits_{l=1}^J(\r_h\tilde{\w}^{(l)}+\nabla \tilde{p}^{(l)}+\tilde{\bf R}^{(l)}_h)\quad\mbox{on}~~G_0.
$$
Then we can build a stable discrete regular decomposition of $\tilde{\v}^{(0)}_h$ on $G_0$ and furthermore build the desired regular decomposition (\ref{regulardecomposition:3.new0}) of $\v_h$ on $G$
as in Lemma \ref{lem:helm5}.
\hfill $\Box$

\section{Proof of Theorem \ref{thm:thm5.1}-Theorem \ref{thm:thm5.4}}
\setcounter{equation}{0}
In this section we are devoted to the proof of Theorem \ref{thm:thm5.1}-Theorem \ref{thm:thm5.4}. We use the same notations as in Section 3 (for example, $\Gamma_i=\Gamma\cap\partial G_i$ for $i=1,2$).
\vskip 0.2in
\no{\it Proof of Theorem \ref{thm:thm5.1}}. We recall that $\bar{G}=\bar{G}_1\cup\bar{G}_2$, and $\bar{G}_1\cap\bar{G}_2$ is the common edge $\E$ of $G_1$ and $G_2$.
We consider three different relations of position between $\E$ and $\Gamma_i$ ($i=1,2$).

(i) $\E\subset\Gamma_i$ for $i=1,2$. In this case, we have $\E=\Gamma_1\cap\Gamma_2$. We use Theorem \ref{teor:helm1} to build a regular decomposition
of $\v_h|_{G_i}$ independently for $i=1,2$. Then the resulting functions $p_{h,i}$ and $\w_{h,i}$ vanish on $\Gamma_i$ and so they also vanish on $\E$. Thus we can directly extend $p_{h,i}$ and $\w_{h,i}$
onto the domain $G\backslash G_i$ by zero to get extension functions $\tilde{p}_{h,i}$ and $\tilde{\w}_{h,i}$ ($i=1,2$). Define $p_h=\tilde{p}_{h,1}+\tilde{p}_{h,2}$ and $\w_h=\tilde{\w}_{h,1}+\tilde{\w}_{h,2}$.
By the two functions, we can build the regular decomposition (\ref{5.new1}) of $\v_h$ on the global $G$. In this case, there is no logarithmic factor in the stability estimates.

(ii) $\E$ is contained in only one of $\Gamma_1$ and $\Gamma_2$, for example, $\E$ is contained in $\Gamma_1$ but is not contained in $\Gamma_2$.
We use Theorem \ref{teor:helm1} to build a regular decomposition of $\v_h|_{G_1}$ associated with $\Gamma_1$, but use Theorem \ref{teor:helm1}
to get a regular decomposition of $\v_h|_{G_2}$ associated with $\Gamma_2\cup\E$ (notice that $\v_h|_{G_2}$ vanishes on $\Gamma_2$ and $\E$ since $\E\subset\Gamma_1$).
It is clear that the resulting functions $p_{h,i}$ and $\w_{h,i}$ vanish on $\Gamma_i$ and $\E$. Then the desired regular decomposition can be built as in the above situation.

(iii) $\E\cap\Gamma_i=\emptyset$ for $i=1,2$. In this case, we have $\Gamma_1\cap\Gamma_2=\emptyset$. If both $\Gamma_1$ and $\Gamma_2$ are non-empty sets, then $\Gamma$ is non-connected;
when $\Gamma$ is connected, one of $\Gamma_1$ and $\Gamma_2$ is connected and the other is empty.

Without loss of generality, we assume that both $\Gamma_1$ and $\Gamma_2$ are non-empty sets. We first use Theorem \ref{teor:helm1} to build a regular decomposition of $\v_h|_{G_1}$
\ee
\v_h=p_{h,1}+\r_h\w_{h,1}+{\bf R}_{h,1}~~~~\mbox{on}~~G_1 \label{5.new3}
\e
with $p_{h,1}\in Z_h(G_1)$ and $\w_{h,1}\in (Z_h(G_1))^3$ vanishing on $\Gamma_1$ (and so ${\bf R}_{h,1}$ has zero degrees of freedom on $\Gamma_1$). Moreover, we have (if $\Gamma_1$ is connected and contains edges)
\ee
\|\w_{h,1}\|_{1,G_1}+h^{-1}\|{\bf R}_{h,1}\|_{0,G_1}\stl\log(1/h)\|\c\, \v_h\|_{0,G_1}
\label{stab:5.new2}\e
and
\ee
\|\w_{h,1}\|_{0,G_1}+\|p_{h,1}\|_{1,G_1}\stl\log(1/h)\|\v_h\|_{\c,G_1}. \label{5.new4}
\e
Notice that $p_{h,1}$, $\w_{h,1}$ and ${\bf R}_{h,1}$ may do not vanish on $\E$ since $\E\cap\Gamma_1=\emptyset$.
We extend $p_{h,1}$, $\w_{h,1}$ and ${\bf R}_{h,1}$ onto $G_2$ such that the resulting extensions $\tilde{p}_{h,1}\in Z_h(G)$, $\tilde{\w}_{h,1}\in (Z_h(G))^3$
and $\tilde{\bf R}_{h,1}\in V_h(G)$ have zero degrees of freedom on the nodes or the edges in $G_2\backslash\E$. Define
\ee
\v^{\ast}_h=\v_h-(\tilde{p}_{h,1}+\r_h\tilde{\w}_{h,1}+\tilde{\bf R}_{h,1})~~~\mbox{on}~~G.\label{5.new5}
\e

Next we build a regular decomposition of $\v^{\ast}_h$. Since $(\tilde{p}_{h,1}+\r_h\tilde{\w}_{h,1}+\tilde{\bf R}_{h,1})\times\n=\v_h\times\n$ on $\partial G_1$,
we have $\lambda_e(\v^{\ast}_h)=0$ for $e\subset\E\subset\partial G_1$.
Moreover, noting that $\Gamma_2\cap\E=\emptyset$, we have $\lambda_e(\v^{\ast}_h)=\lambda_e(\v_h)$ for $e\subset\Gamma_2\subset\partial G_2\backslash\E$, which implies that
$\lambda_e(\v^{\ast}_h)=0$ for $e\subset\Gamma_2$. In summary, $\v_h^{\ast}$ has zero degrees of freedom on $\E\cup\Gamma_2$. Thus, by Theorem \ref{teor:helm1} the function $\v^{\ast}_h$ admits a regular decomposition
\ee
\v^{\ast}_h=p^{\ast}_{h,2}+\r_h\w^{\ast}_{h,2}+{\bf R}^{\ast}_{h,2}~~~\mbox{on}~~G_2, \label{5.new6}
\e
with $p^{\ast}_{h,2}\in Z_h(G_2)$ and $\w^{\ast}_{h,2}\in (Z_h(G_2))^3$ vanishing on $\E\cup\Gamma_2$. Moreover, we have
(the set $\E\cup\Gamma_2$ is non-connected since $\E\cap\Gamma_2=\emptyset$ and $\Gamma_2\not=\emptyset$)
\ee
\|\w^{\ast}_{h,2}\|_{1,G_2}+h^{-1}\|{\bf R}^{\ast}_{h,2}\stl\log(1/h)\|\v^{\ast}_h\|_{\c,G_2}
\label{stab:5.new3}\e
and
\ee
\|\w^{\ast}_{h,2}\|_{0,G_2}+\|p^{\ast}_{h,2}\|_{1,G_2}\stl\log(1/h)\|\v^{\ast}_h\|_{\c,G_2}. \label{5.new7}
\e
When $\Gamma$ is connected, for example, $\Gamma_1$ is connected but $\Gamma_2=\emptyset$, the norm on the right side of (\ref{stab:5.new3}) can be replaced by the $\c$ semi-norm.

Since $p^{\ast}_{h,2}$, $\w^{\ast}_{h,2}$ and ${\bf R}^{\ast}_{h,2}$ have zero degrees of freedom on $\E$, these functions can be extended onto $G_1$ by zero and the resulting
extensions $\tilde{p}^{\ast}_{h,2}$, $\tilde{\w}^{\ast}_{h,2}$ and $\tilde{\bf R}^{\ast}_{h,2}$ satisfy $\tilde{p}^{\ast}_{h,2}\in Z_h(G)$, $\tilde{\w}^{\ast}_{h,2}\in (Z_h(G))^3$
and $\tilde{\bf R}^{\ast}_{h,2}\in V_h(G)$, respectively. Define
$$ p_h=\tilde{p}_{h,1}+\tilde{p}^{\ast}_{h,2},~\w_h=\tilde{\w}_{h,1}+\tilde{\w}^{\ast}_{h,2}~~\mbox{and}~~ {\bf R}_h=\tilde{\bf R}_{h,1}+\tilde{\bf R}^{\ast}_{h,2}.$$
Then $p_h$, $\w_h$ and ${\bf R}_h$ have zero degrees of freedom on $\Gamma=\Gamma_1\cup\Gamma_2$ and the regular decomposition (\ref{5.new1}) follows by (\ref{5.new5}) and (\ref{5.new6}).

In an analogous way with Step 2 in the proof of Theorem 4.1, we can verify the estimates (\ref{stab:5.new1}) and (\ref{5.new2}) by using (\ref{stab:5.new2})-(\ref{5.new4}) and (\ref{stab:5.new3})-(\ref{5.new7}).
\hfill $\Box$

In the following we want to prove Theorem \ref{thm:thm5.2} (and Theorem 3.4). For this purpose, we first prove two auxiliary results involving one vertex on a usual polyhedron domain.

\begin{lemma}\label{lem:splitV} Let $G$ be a usual polyhedron and $\vv$ be a vertex of $G$. Assume that $\v_h$ is a function in $V_h(G)$. Then we can write $\v_h$ as
$$ \v_h=\nabla p_h+\r_h\w_h+{\bf R}_h$$
for some $p_h\in Z_h(G)$, $\w_h\in (Z_h(G))^3$ and $
{\bf R}_h\in V_h(G)$ satisfying $p_h(\vv)=0$ and $\w_h(\vv)=\0$. Moreover, we have
\ee
\|\w_h\|_{1,G}+h^{-1}\|{\bf R}_h\|_{0,G}\stl\log(1/h)\|\c\,\v_h\|_{0,G}\label{stab:4.new16}\e
and
\ee
\|\w_h\|_{0,G}+\|\nabla p_h\|_{0,G}\stl\log(1/h)\|\v_h\|_{\c,G}. \label{stab:4.new17}\e
\end{lemma}

\no {\it Proof}. Consider a (closed) face $\ff$ containing $\vv$ as its vertex.
Like Step 1 in the proof of Lemma~\ref{lem:helm3}, we can define
$\phi_{\partial\ff}$ to be a function that is piecewise linear and
continuous on $\partial\ff$ such that $\phi_{\partial\ff}(\vv)=0$,
and define $C_{\partial\ff}$ to be a constant such that
$\v_h\cdot\t_{\partial\ff}=\phi'_{\partial\ff}+C_{\partial\ff}$ on
$\partial\ff$. In fact, they can be defined as
\begin{equation}
C_{\partial\ff} = \frac{1}{l}\int_0^{l}(\v_h\cdot\t_{\partial\ff})(s)\,ds, \quad
\phi_{\partial\ff}(t)
=\int_0^{t} (\v_h\cdot\t_{\partial\ff}-C_{\partial\ff})(s)ds\,,
\quad \forall\,t\in [0, l]\,, \label{6.new00}
\end{equation}
where $l$ is the length of $\partial\ff$ and $t=0$ (and $t=l$) corresponds to the vertex $\vv$.
Let $c=\gamma_{\E}(\phi_{\partial\ff})$ denote the average of $\phi_{\partial\ff}$ on $\E$, where $\E$ is
an edge or a union of several edges of $\ff$. Define an extension $\tilde{\phi}_{\vv}\in Z_h(G)$ of $\phi_{\partial\ff}$, such that
$\tilde{\phi}_{\vv}$ equals to the average $c=\gamma_{\E}(\phi_{\partial\ff})$ at all of the nodes
on $G$ except those on $\partial\ff$. Then (cf.\,Lemma 7.6 and Lemma 7.1 of \cite{Tos})
\ee \|\nabla \tilde{\phi}_{\vv}\|_{0,G}= \|\nabla(\tilde{\phi}_{\vv}-c)\|_{0,G}\stl\|\phi_{\partial\ff}-c\|_{0,\partial\ff}\stl\|\phi'_{\partial\ff}\|_{H^{-1}(\partial\ff)}
\stl\log(1/h)\|\v_h\|_{\c,~G}.\label{stab:4.new017}\e
We define a similar extension $\tilde C_{\vv}$ of $C_{\partial\ff}$ with $\tilde{C}_{\E}$ (defined in the proof of Lemma~\ref{lem:helm3}), such that
$\tilde C_{\vv}$ belongs to $(Z_h(G))^3$ and vanishes at
$\vv$, and it satisfies the condition $(\r_h\tilde
C_{\vv})\cdot\t_{\partial\ff}=C_{\partial\ff}$ on $\partial\ff$ and the stability
$$ \|\tilde C_{\vv}\|_{1,G}\stl\|\tilde C_{\vv}\|_{0,\partial\ff}\stl |C_{\partial\ff}|\stl\log^{{1\over 2}}(1/h)\|\c~\v_h\|_{0,G}. $$
Define
$$ \hat{\v}_{h,\vv}=\v_h-(\nabla \tilde{\phi}_{\vv}+\r_h\tilde C_{\vv}). $$
Then we have $\hat{\v}_{h,\vv}\cdot\t_{\partial\ff}=0$. As in the proof of Lemma~\ref{lem:helm3}, we can use Lemma \ref{lem:helm2} for
$\hat{\v}_{h,\vv}$ to build the desired decomposition of $\v_h$. \hfill $\Box$

\begin{remark} Differently from the second inequality in (\ref{stab:4.new14}), the inequality (\ref{stab:4.new17})
holds only for the $H^1$ semi-norm of $p_h$. The main reason is that a stable estimate of $\|\tilde{\phi}_{\vv}\|_{0,G}$ cannot be built except
that the constant $c$ in (\ref{stab:4.new017}) vanishes. In the above proof, the face $\ff$ can be replaced by a connected ``Lipschitz" union of some faces of $G$,
provided that union set has one connected boundary containing $\vv$ as one of its vertices (the union cannot be $\partial G$ itself).
\end{remark}
{\bf Corollary 6.1}. Let $\bar{G}=\bar{G}_1\cup\bar{G}_2$ with $G_i$ be a usual polyhedron domain. Assume that $\bar{G}_1\cap\bar{G}_2$ is just the common vertex $\vv$ of $G_1$ and $G_2$.
Then any function $\v_h\in V_h(G)$ admits a discrete regular decomposition on $G$ and the resulting functions have the stability estimates given in Lemma \ref{lem:splitV}.\\
\no{\it Proof}. We use Lemma \ref{lem:splitV} to build a regular decomposition of $\v_h|_{G_i}$ independently for $i=1,2$. Since the resulting nodal finite element functions vanish at $\vv$,
the two regular decompositions naturally define a regular decomposition for $\v_h$ on the global domain $G$.   \hfill $\Box$


In the rest of this section, for a face (or a connected ``Lipschitz" union of some faces) $\ff$, we always use $\phi_{\partial\ff}$ to denote the function defined in the proof of Lemma \ref{lem:splitV}.

\begin{lemma} \label{lem:helmFV} Let $G$ be a usual polyhedron, and let $\Gamma$ be a (closed) union of some faces and edges of $G$, and $\vv$ be a vertex of $G$ ($\vv\not\in\Gamma$).
Assume that $\v_h$ satisfies {\bf Assumption 3.1}. Suppose that there exists a face (or a connected ``Lipschitz" union of some faces) $\ff$ containing $\vv$ such that $\gamma_{\E}(\phi_{\partial\ff})=0$
for an edge $\E$ of $\ff$, where either $\ff\cap\Gamma=\emptyset$ or $\partial\ff\cap\Gamma$ is a union of edges of $G$ (in this case, we require that $\E\subset\partial\ff\cap\Gamma$). Then $\v_h$ can be decomposed as
\ee \v_h=\nabla p_h+\r_h\w_h+{\bf R}_h\label{5.new0007}
\e
for some $p_h\in Z_h(G)$, $\w_h\in (Z_h(G))^3$ and
${\bf R}_h\in V_h(G)$ such that $p_h$ and $\w_h$ vanish on $\Gamma$ and $\vv$. Moreover, we have
 \ee \|\w_h\|_{1,G}+h^{-1}\|{\bf R}_h\|_{0,G}\stl\log(1/h)\|\v_h\|_{\c,G}
\label{5.new0008}\e
and
\ee
\|\w_h\|_{0,G}+\|p_h\|_{1,G}\stl\log(1/h)\|
\v_h\|_{\c,G}. \l{5.new0009}\e
In particular, when $\Gamma$ is connected and $\E\subset\partial\ff\cap\Gamma$, the norm on the right side of (\ref{5.new0008}) can be replaced by the $\c$ semi-norm.
\end{lemma}
\no {\it Proof}. A key step in the proof is to define two functions $\tilde{\phi}_{\vv}\in Z_h(G)$ and $\tilde C_{\vv}\in (Z_h(G))^3$ such that they have zero degrees of freedom
at $\vv$ and $\Gamma$ and satisfy $\v_h\cdot\t_{\partial\ff}=(\nabla\tilde{\phi}_{\vv}+\r_h\tilde C_{\vv})\cdot\t_{\partial\ff}$ on $\partial\ff$.

We first consider the case with $\ff\cap\Gamma=\emptyset$. For this case, we can simply define $\tilde{\phi}_{\vv}$ and $\tilde C_{\vv}$ as in the proof of Lemma \ref{lem:splitV}. Set
$$ \hat{\v}_{h,\vv}=\v_h-(\nabla\tilde{\phi}_{\vv}+\r_h\tilde C_{\vv}). $$
Then we have $\lambda_e(\hat{\v}_{h,\vv})=0$ for any $e\subset\partial\ff\cup\Gamma$ by the definitions of $\tilde{\phi}_{\vv}$ and $\tilde C_{\vv}$, together with the assumption on $\v_h$.
Thus we can use Theorem \ref{teor:helm1} (with $\Gamma$ replacing by $\partial\ff\cup\Gamma$) to build a regular decomposition of $\hat{\v}_{h,\vv}$ such that the functions defined
by the decomposition have zero degrees of freedom on $\partial\ff\cup\Gamma$, which implies that the two $H^1$ finite element functions vanish at $\vv$ (and $\Gamma$) since $\vv\in\partial\ff$.
Based on this decomposition we further get the desired decomposition of $\v_h$ and the estimates (refer to the proof of Lemma~\ref{lem:splitV}).
Here we need to use the assumption $\gamma_{\E}(\phi_{\partial\ff})=0$ to get the $L^2$ stability of $\tilde{\phi}_{\vv}$ (see Remark 6.1).

Now we investigate the case with $\partial\ff\cap\Gamma$ being a union of edges of $G$. For this case, the analysis is a bit complicated
because the functions $\tilde{\phi}_{\vv}$ and $\tilde{C}_{\vv}$ defined in the proof of Lemma \ref{lem:splitV} may do not vanish on $\partial\ff\cap\Gamma$ although
they vanish at $\vv$. Based on this consideration, we have to define two new functions $\tilde{\phi}_{\vv}$ and $\tilde{C}_{\vv}$ such that they can satisfy the requirements mentioned above.

\begin{center}
\includegraphics[width=6cm,height=6cm]{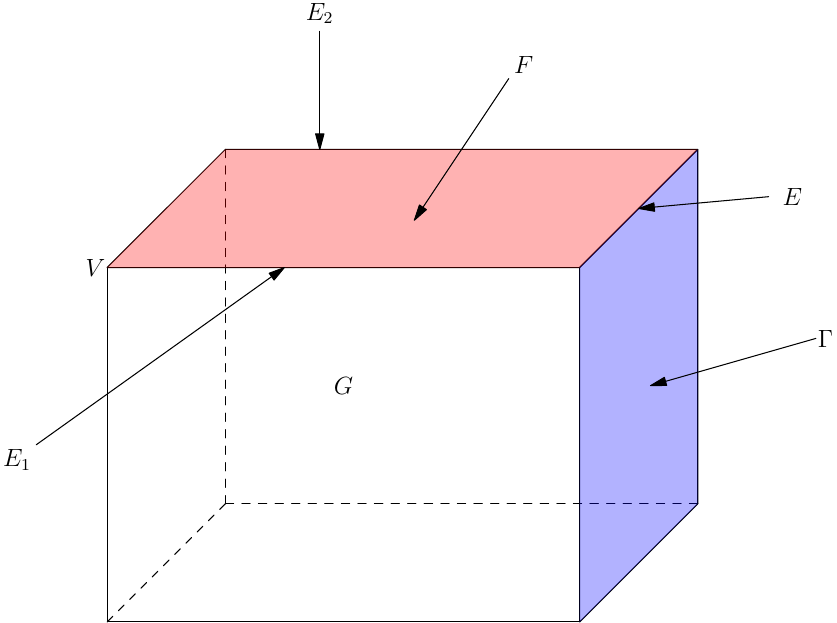}
\vskip 0.1in
\centerline{Figure 8: An example of cube: the blue face denote $\Gamma$, the red face denotes $\ff$, $\ff\cap\Gamma=\E$.}
\end{center}

Without loss of generality, we assume that $\ff$ is a face and $\partial\ff\cap\Gamma=\E$ is an edge of $G$ (see Figure 8).
For more general situation, there is no essential change in the analysis. Let $\phi_{\partial\ff}$ and $C_{\partial \ff}$ be defined by (\ref{6.new00}).
 The basic ideas are to choose a piecewise constant function $\varepsilon_{\partial\ff}$ on $\partial\ff$ and define
$\tilde{\phi}_{\vv}\in Z_h(G)$ such that $\tilde{\phi}_{\vv}$ vanishes at all of the nodes on $\bar{G}\backslash\partial\ff$ and is given by the following formula on $\partial\ff$
\ee
\tilde{\phi}_{\vv}(t)=\phi_{\partial\ff}(t)-\int_0^t\varepsilon_{\partial\ff}ds,~~~~t\in[0,|\partial\ff|], \label{6.new000}
\e
where $t=0$ and $t=|\partial\ff|$ are the arc-length coordinates of the vertex $\vv$.

In the following we give a definition of $\varepsilon_{\partial\ff}$.
Notice that $\phi_{\partial\ff}(|\partial\ff|)=0$, the function $\varepsilon_{\partial\ff}$ should satisfy $\int_0^{|\partial\ff|}\varepsilon_{\partial\ff}ds=0$ so that
$\tilde{\phi}_{\vv}(|\partial\ff|)=0$. Besides, since we hope that $\tilde{\phi}_{\vv}$ vanishes on $\E$, we naturally define $\varepsilon_{\partial\ff}=\phi'_{\partial\ff}$ on $\E$.
Let $\E_1,~\E_2\subset\partial\ff$ be the two neighboring edges with $\E$ (see Figure 8), and let $t_1$ and $t_2$ denote the arc-length coordinates of the two endpoints of $\E$ ($t_1<t_2$).
Then the piecewise constant function $\varepsilon_{\partial\ff}$ can be defined as
\[\varepsilon_{\partial\ff}=\cases{\quad\quad\phi'_{\partial\ff}\quad\quad\mbox{on}\quad \E,\cr\phi_{\partial\ff}(t_1)/|\E_1|\quad\mbox{on}\quad \E_1,\cr
-\phi_{\partial\ff}(t_2)/|\E_2|\quad\mbox{on}\quad \E_2,\cr\quad\quad 0\quad\mbox{on}\quad \partial\ff\backslash(\E\cup\E_1\cup\E_2).}\]
It can be verified that, for the above $\varepsilon_{\partial\ff}$, the function $\tilde{\phi}_{\vv}$ defined in (\ref{6.new000}) indeed vanishes on $\vv$ and $\E$.

For the new function $\tilde{\phi}_{\vv}$, we need to define the corresponding function $\tilde{C}_{\vv}$, which can be regarded as a variant of $C_{\partial\ff}$.
Since $\E\subset\partial\ff\cap\Gamma$ and $\v_h$ has zero degrees of freedom on $\Gamma$, we have $(\v_h\cdot\t_{\partial\ff})|_{\E}=\v_h\cdot\t_{\E}=0$. Then, by the direct calculation, we get
$$ \phi_{\partial\ff}(t)=\int_0^{t}(\v_h\cdot\t_{\partial\ff}-C_{\partial\ff})ds=\int_0^{t_1}\v_h\cdot\t_{\partial\ff}ds-t C_{\partial\ff},~~~~t\in[t_1, t_2]. $$
Thus we have $\varepsilon_{\partial\ff}=\phi'_{\partial\ff}=-C_{\partial\ff}$ on $\E$. Define $\tilde{C}_{\vv}\in(Z_h(G))^3$ such that $\tilde{C}_{\vv}\cdot\t_{\partial\ff}=C_{\partial\ff}+\varepsilon_{\partial\ff}$ on $\partial\ff$, and $\tilde{C}_{\vv}$ vanishes at $\vv$ and all of the nodes on $\bar{G}\backslash\partial\ff$. It is clear that $\tilde{C}_{\vv}\cdot\t_{\partial\ff}=0$ on $\E$ by the relation $\varepsilon_{\partial\ff}=-C_{\partial\ff}$. Moreover, the following decomposition holds for the modified functions $\tilde{\phi}_{\vv}$ and $\tilde{C}_{\vv}$
\ee
\v_h\cdot\t_{\partial\ff}=C_{\partial\ff}+\phi'_{\partial\ff}=\tilde{C}_{\vv}\cdot\t_{\partial\ff}+\tilde{\phi}'_{\vv}~~~\mbox{on}~~\partial\ff. \label{equality-5.0002}
\e

From the above discussions, we know that the functions $\tilde{C}_{\vv}$ and $\tilde{\phi}_{\vv}$ have zero degrees of freedom on $\Gamma$ and $\vv$
(since the two functions naturally have zero degrees of freedom on $\Gamma\backslash\E$ by their definitions).
In addition, the functions $\tilde{C}_{\vv}$ and $\tilde{\phi}_{\vv}$
possess the same stabilities as the corresponding functions defined in the proof of Lemma \ref{lem:splitV}. Here we use the condition that $\tilde{\phi}_{\vv}$ vanishes on
 $\E$
to get the $L^2$ stability of $\tilde{\phi}_{\vv}$.
Define
$$ \hat{\v}_{h,\vv}=\v_h-(\nabla\tilde{\phi}_{\vv}+\r_h\tilde C_{\vv}). $$
By (\ref{equality-5.0002}) and the assumption on $\v_h$, we have $\lambda_e(\hat{\v}_{h,\vv})=0$ for any $e\subset\partial\ff\cup\Gamma$. Thus we can prove the desired results as in the case with $\ff\cap\Gamma=\emptyset$ (see the previous part). In particular, if $\Gamma$ is connected, then $\ff\cup\Gamma$
is also connected, and so the complete norm on the right side of (\ref{5.new0008}) can be replaced by the $\c$ semi-norm (see the conclusion given at the bottom of Theorem 3.1).
\hfill$\Box$


\begin{remark} We find, from Lemma \ref{lem:helmFV}, that a vertex is essentially different from an edge. This phenomenon can be intuitively explained as follows: the value of an edge finite element
function $\v_h$ is not uniquely defined at a vertex, but the two nodal finite element functions defined by the regular decomposition of $\v_h$ are required to vanish at the vertex. Thus there is a gap between $\v_h$
and the gradient of the scalar finite element function, which needs to be filled by a constraint of $\v_h$.
\end{remark}


Now we prove Theorem \ref{thm:thm5.2} by Theorem \ref{teor:helm1}, Lemma \ref{lem:splitV} and Lemma \ref{lem:helmFV}.
\vskip 0.2in
\no{\it Proof of Theorem \ref{thm:thm5.2}}. We recall that $\bar{G}=\bar{G}_1\cup\bar{G}_2$, where $G_1$ and $G_2$ are
two usual polyhedra that intersect at a vertex $\vv$, and $\Gamma_i=\Gamma\cap\partial G_i$. We consider three different situations, respectively.

(i) Both $\Gamma_1$ and $\Gamma_2$ contain the vertex $\vv$, which implies that the intersection of $\Gamma_1$ and $\Gamma_2$ is just $\vv$.

We use Theorem \ref{teor:helm1} to build a regular decomposition of $\v_h|_{G_i}$ independently for $i=1,2$, such that
the resulting nodal finite element functions have zero degrees of freedom on $\Gamma_i$ ($i=1,2$).
Since $\vv=\Gamma_1\cap\Gamma_2$, the nodal finite element functions also vanish at $\vv$.
Thus the two regular decompositions can naturally define the desired regular decomposition on the global domain $G$. In this case,
the constraint is unnecessary, i.e., we can choose ${\mathcal F}$ as the zero functional.

(ii) One of $\Gamma_1$ and $\Gamma_2$ contains the vertex $\vv$, for example, $\vv\in\Gamma_1$ but $\vv\notin\Gamma_2$.

We choose a face $\ff$ of $G_2$ such that $\vv\in\ff$ (see the left graph in Figure 9). Assume that either $\ff\cap\Gamma_2$ is an empty set or $\ff\cap\Gamma_2$ is a union of edges.
Let $s$ be the arc-length variable along $\partial\ff$, taking values from
$0$ to $|\partial\ff|$, which is the total length of $\partial\ff$. Define
$$
C_{\partial\ff} = \frac{1}{|\partial\ff|}\int_{\partial\ff}(\v_h\cdot\t_{\partial\ff})(s)\,ds, \quad
\phi_{\partial\ff}(t)
=\int_0^{t} (\v_h\cdot\t_{\partial\ff}-C_{\partial\ff})(s)ds\,,
\quad \forall\,t\in [0, |\partial\ff|].
$$

For an edge $\E$ of $\ff$, let $\gamma_{\E}(\phi_{\partial\ff})$ denote the integration average
of $\phi_{\partial\ff}$ on $\E$. Define the functional ${\cal F}$ as ${\cal F}\v_h=\gamma_{\E}(\phi_{\partial\ff})$, where
we require that $\E\subset\ff\cap\Gamma_2$ in the case of $\ff\cap\Gamma_2\not=\emptyset$. Then the constraint ${\cal F}\v_h=0$ means that
$\gamma_{\E}(\phi_{\partial\ff})=0$.

\begin{center}
\includegraphics[width=6.3cm,height=6cm]{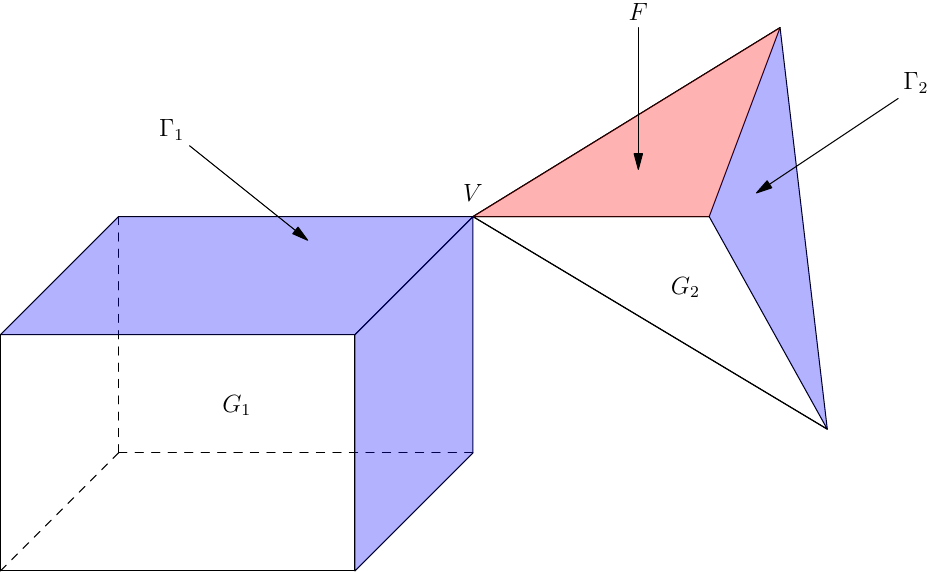}\quad\quad\includegraphics[width=6.3cm,height=6cm]{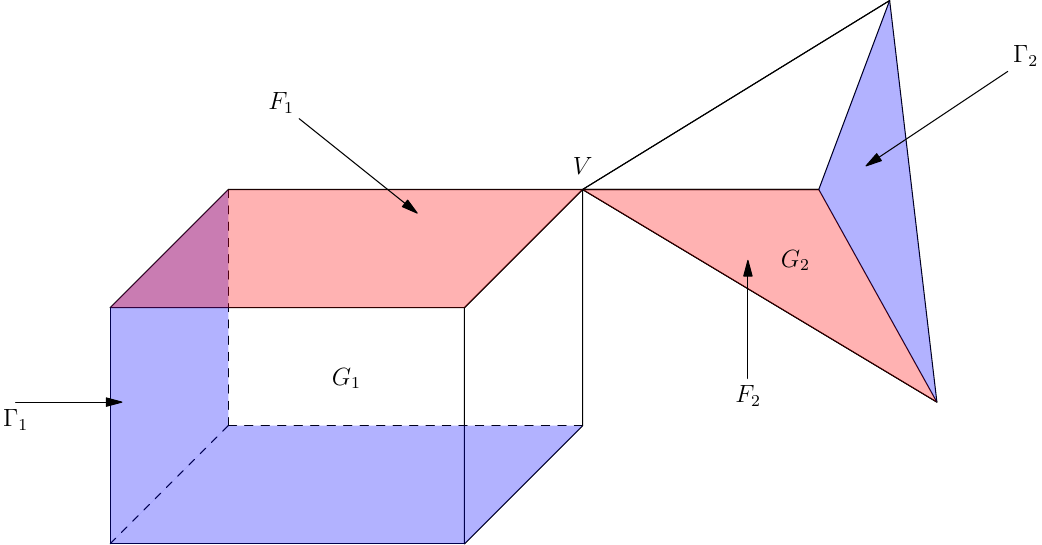}
\end{center}
Figure 9. The blue part of $\partial G_i$ denotes $\Gamma_i$. Situation (i) is shown by the left graph: $\vv\in\Gamma_1$, $\vv\notin\Gamma_2$, $\ff\subset\partial G_2$ and $\ff\cap\Gamma_2$ is an edge,
the red part of $\partial G_2$ denotes $\ff$;
Situation (ii) is shown by the right graph: $\vv\notin\Gamma_1\cup\Gamma_2$, $\ff_1\subset\partial G_1$ and $\ff_2\subset\partial G_2$, the red part of $\partial G_i$ denotes $\ff_i$.

\vskip 0.2in
We use Theorem 3.1 to build a discrete regular decomposition on $G_1$ such that
the resulting functions have zero degrees of freedom on $\Gamma_1$. Since $\vv\in\Gamma_1$, the $H^1$ finite element functions vanish at $\vv$. Moreover, we use Lemma \ref{lem:splitV}
(if $\Gamma_2=\emptyset$) or Lemma \ref{lem:helmFV} (if $\Gamma_2\not=\emptyset$) to build a discrete regular decomposition on $G_2$ such that
the resulting three functions have zero degrees of freedom on $\Gamma_2$ and the two $H^1$ finite element functions vanish at $\vv$. Then the two regular decompositions
naturally define the desired regular decomposition on $G$. The stability of the $L^2$ norm of $p_h$ can be guaranteed by the condition $\gamma_{\E}(\phi_{\partial\ff})=0$.

(iii) Neither $\Gamma_1$ nor $\Gamma_2$ contains the vertex $\vv$, i.e., $\vv\notin\Gamma=\Gamma_1\cup\Gamma_2$.

We choose two faces $\ff_1\subset\partial G_1$ and $\ff_2\subset\partial G_2$ such that both of them
contain $\vv$ (see the right graph in Figure 9). Suppose that either $\ff_i\cap\Gamma_i$ is an empty set or $\partial\ff_i\cap\Gamma_i$ is a union of some edges of $G_i$ ($i=1,2$).
The face $\ff_i$ can be replaced by a connected ``Lipschitz" union of some faces of $G_i$,
provided that the union has one connected boundary containing $\vv$ as one of its vertices (the union cannot be $\partial G_i$ itself).
Define
$$
C_{\partial\ff_i} = \frac{1}{|\partial\ff_i|}\int_{\partial\ff_i}(\v_h\cdot\t_{\partial\ff_i})(s)\,ds, \quad
\phi_{\partial\ff_i}(t)
=\int_0^{t} (\v_h\cdot\t_{\partial\ff_i}-C_{\partial\ff_i})(s)ds+c_i\,,
\quad \forall\,t\in [0, |\partial\ff_i|]\,,
$$
where the constant $c_i$ is chosen
such that $\gamma_{\E_i}(\phi_{\partial\ff_i})=0$ for some edge $\E_i\subset\partial\ff_i$. When $\partial\ff_i\cap\Gamma_i$ is a union of some edges of $G_i$,
we require that $\E_i\subset\partial\ff_i\cap\Gamma_i$.

Define the functional ${\cal F}$ as ${\cal F}\v_h=\phi_{\partial\ff_1}(|\partial\ff_1|)-\phi_{\partial\ff_2}(|\partial\ff_2|)$. The constraint ${\cal F}\v_h=0$ means
that $\phi_{\partial\ff_1}(|\partial\ff_1|)=\phi_{\partial\ff_2}(|\partial\ff_2|)$.

For convenience, we only consider the case with $\ff_i\cap\Gamma_i=\emptyset$ ($i=1,2$), otherwise, we need to repeat the technique in Lemma \ref{lem:helmFV}. As in the proof of Lemma \ref{lem:splitV},
we respectively define the extensions $\tilde{\phi}^{(i)}_{\vv}$ and
$\tilde{C}^{(i)}_{\vv}$ of $\phi_{\partial\ff_i}$ and $C_{\partial\ff_i}$, but define the values of $\tilde{\phi}^{(i)}_{\vv}$ as zero at all of the nodes in $G_i$ except on $\partial\ff_i$.
Set
$$ \hat{\v}^{(i)}_{h,\vv}=\v_h|_{G_i}-(\tilde{\phi}^{(i)}_{\vv}+\r_h\tilde{C}^{(i)}_{\vv})\quad\mbox{on}~~G_i. $$
Then we have $\lambda_e(\hat{\v}^{(i)}_{h,\vv})=0$ for any $e\subset\partial\ff_i\cup\Gamma_i$ by the definitions of $\tilde{\phi}^{(i)}_{\vv}$ and $\tilde{C}^{(i)}_{\vv}$.
Thus we can use Theorem 3.1 for $\hat{\v}^{(i)}_{h,\vv}$ to build a regular decomposition of $\v_h|_{G_i}$ and the resulting functions have the corresponding stability estimates by the condition $\gamma_{\E_i}(\phi_{\partial\ff_i})=0$
(refer to the proof of Lemma~\ref{lem:splitV} and Remark 6.1). Using the definitions of the functions $p_{h,i}$ and $\w_{h,i}$ in the regular decomposition of $\v_h|_{G_i}$, together with the condition $\phi_{\partial\ff_1}(|\partial\ff_1|)=\phi_{\partial\ff_2}(|\partial\ff_2|)$, yields $p_{h,1}=p_{h,2}$ and $\w_{h,i}=\0$ at the vertex $\vv$. Therefore, we can naturally define a regular decomposition of $\v_h$ on
the global $G$ and obtain the desired stability estimates.
\hfill $\Box$
\vskip 0.2in
\no{\it Proof of Theorem \ref{thm:thm5.4}}. Let $\ff_i$ be a face (or a connected ``Lipschitz" union of faces) satisfying $\vv\in\partial\ff_i\subset\partial G_i$
and define a function $\phi_{\partial\ff_i}$ as in the proof of Theorem \ref{thm:thm5.2} ($i=1,\cdots,s$). We define the functional ${\cal F}_i$ by
$$ {\cal F}_i\v_h=\phi_{\partial\ff_{i+1}}(\vv)-\phi_{\partial\ff_i}(\vv)\quad (i=1,\cdots,s-1);$$
Then we can prove the desired results in an analogous way with the proof of Theorem \ref{thm:thm5.2}.

\hfill $\Box$\\
{\bf Acknowledgments}. The author would like to thank the anonymous reviewer, who gives many insightful comments to improve the presentation of this paper.
Moreover, the author wishes to thank Professor Jinchao Xu for suggesting him to study the convergence of the HX preconditioner for the case with jump coefficients.
In particular, Professor Jinchao Xu suggested the author to consider regular decomposition on non-Lipchitz domains, which makes the results in this paper be richer and more complete.

\end{document}